\documentclass[11pt]{article}

\usepackage{amsmath, graphicx, latexsym, amssymb, amsthm, amsfonts}
\usepackage[mathscr]{eucal}
\usepackage{graphics}
\usepackage{epsfig}

\DeclareMathOperator*{\argmin}{arg\,min}

\newtheorem{lemma}{Lemma}
\newtheorem{theorem}{Theorem}
\newtheorem{corollary}{Corollary}
\newtheorem{definition}{Definition}
\newtheorem{example}{Example}
\newtheorem{proposition}{Proposition}

\newtheorem{conjecture}{Conjecture}

\newtheorem{remark}{Remark}

\newcommand{\beginsec}{
\setcounter{lemma}{0}
\setcounter{theorem}{0}
\setcounter{corollary}{0}
\setcounter{definition}{0}
\setcounter{example}{0}
\setcounter{proposition}{0}
\setcounter{condition}{0}
\setcounter{assumption}{0}
\setcounter{conjecture}{0}
\setcounter{problem}{0}
\setcounter{remark}{0}
}

\newcommand{\la}{\lambda}
\newcommand{\eps}{\varepsilon}
\newcommand{\ph}{\varphi}

\newcommand{\al}{\alpha}

\newcommand{\gam}{\gamma}

\newcommand{\sig}{\sigma}
\newcommand{\del}{\delta}

\newcommand{\Gam}{\mathnormal{\Gamma}}
\newcommand{\Del}{\mathnormal{\Delta}}
\newcommand{\Th}{\mathnormal{\Theta}}

\newcommand{\X}{\mathnormal{\Xi}}

\newcommand{\Om}{\mathnormal{\Omega}}

\newcommand{\D}{{\mathbb D}}

\newcommand{\N}{{\mathbb N}}

\newcommand{\R}{{\mathbb R}}

\newcommand{\EE}{{\mathbb E}}
\newcommand{\E}{{\mathbb E}}

\newcommand{\PP}{{\mathbb P}}

\newcommand{\calA}{{\cal A}}
\newcommand{\calB}{{\cal B}}

\newcommand{\calD}{{\cal D}}

\newcommand{\calF}{{\cal F}}
\newcommand{\calG}{{\cal G}}
\newcommand{\calH}{{\cal H}}
\newcommand{\calI}{{\cal I}}

\newcommand{\calL}{{\cal L}}

\newcommand{\calU}{{\cal U}}

\newcommand{\calX}{{\cal X}}

\newcommand{\bB}{{\mathbf B}}

\newcommand{\bT}{{\mathbf T}}

\newcommand{\bx}{{\mathbf x}}

\newcommand{\scrI}{\mathscr{I}}

\renewcommand{\proof}{\noindent{\bf Proof.\ }}

\newcommand{\uu}{\underline}
\newcommand{\oo}{\overline}
\newcommand{\skp}{\vspace{\baselineskip}}

\newcommand{\diag}{{\rm diag}}

\newcommand{\w}{\wedge}

\newcommand{\pl}{\partial}

\newcommand{\To}{\Rightarrow}
\newcommand{\dist}{{\rm dist}}

\newcommand{\iy}{\infty}

\newcommand{\be}{\begin{equation}}
\newcommand{\ee}{\end{equation}}

\newcommand{\noi}{\noindent}
\newcommand{\ds}{\displaystyle}

\begin{document}

\title{An asymptotic optimality result for the multiclass queue with finite buffers
in heavy traffic\thanks{Research
supported in part by the ISF (Grant 1349/08), the US-Israel BSF (Grant 2008466),
and the Technion fund for promotion of research}}
\author{Rami Atar \hspace{4em} Mark Shifrin\\ \\ \\
Department of Electrical Engineering\\
Technion--Israel Institute of Technology\\
Haifa 32000, Israel}

\date{August 21, 2013; revised May 5, 2014}

\maketitle

\begin{abstract}
For a multiclass G/G/1 queue with finite buffers,
admission and scheduling control, and holding and
rejection costs, we construct a policy that is asymptotically optimal in the heavy traffic limit.
The policy is specified in terms of a single parameter which constitutes the free boundary
point from the Harrison-Taksar free boundary problem,
but otherwise depends ``explicitly'' on the problem data.
The {\it $c\mu$ priority rule} is also used by the policy, but in a way that is novel, and, in particular,
different than that used in problems with infinite buffers.
We also address an analogous problem where buffer constraints
are replaced by throughput time constraints.

\skp

\noi{\bf AMS subject classifications:} 60F17, 60J60, 60K25, 93E20

\skp

\noi{\bf Keywords:}
Multiclass G/G/1 queue; Brownian control problems; Bellman equation;
The Harrison-Taksar free boundary problem;
State dependent priorities
\end{abstract}

\section{Introduction}
\beginsec

In this work we consider the problem of finding asymptotically
optimal (AO) controls for the multiclass G/G/1 queue with finite buffers, in heavy traffic.
Upon arrival of a class-$i$ customer into queue $i$ (with
$i\in\{1,\ldots,I\}$ and where $I$ denotes the number of classes),
a decision maker may either accept or reject the job.
In addition, the decision maker controls
the fraction of effort devoted by the server to the customer at the head of queue $i$,
for each $i$.
We refer to the two elements of the control as {\it admission control} and {\it scheduling},
respectively.
The problem considered is to minimize a combination of holding and rejection costs.
The term `heavy traffic' refers to assuming
a critical load condition and observing the model at diffusion scale.
Our interest in this problem stems from recent developments in
the application area of cloud computing. In a hybrid cloud
where a private cloud (namely, a local server) has a given capacity and
memory limits, tasks that cannot be queued in real time are rejected from
the local system and sent to a public cloud, where a fixed charge per usage applies.
For further details on modeling toward these applications,
see \cite{shi-ata-cid}. For a more general modeling framework of data centers, see \cite{ben-men}.
The analysis of the model leads in the scaling limit
to a control problem associated with Brownian motion (BM), often referred to
in this context as a {\it Brownian control problem} (BCP).
Our main result is the convergence of the {\it queueing control problem} (QCP) value function
to that of the BCP, and the construction of a particular AO admission/scheduling policy.
The policy is specified in terms of a free boundary point that is used in solving the BCP,
but otherwise it depends explicitly on the problem data.

A line of research starting from Harrison \cite{har1988} and continuing with
Harrison and van Mieghem \cite{Har-Van},
Harrison \cite{har2000}, \cite{har2003} and Harrison and Williams \cite{Har-Wil} has treated
BCP associated with a broad family of models called
{\it stochastic processing networks}.
These problems, aimed at describing the heavy traffic limits of QCP, were
shown to be equivalent to {\it reduced BCP} (RBCP), in which {\it workload}
plays the role of a state process. RBCP simplify BCP in two ways: Their
state lies in lower dimension, and their form, specifically, that of a singularly controlled
diffusion, makes control theoretic tools applicable.
Addressing these models at the same level of generality,
Atar and Budhiraja \cite{AB2006} and Atar, Budhiraja and Williams \cite{ABW}
use such control theoretic tools
to characterize the BCP (equivalently, RBCP) value functions as solutions
to {\it Hamilton-Jacobi-Bellman} (HJB) equations,
and Budhiraja and Ghosh \cite{BG2006} and \cite{BG2012} prove convergence of QCP value
function to BCP value functions.
Many other works address these models in situations where the BCP are explicitly solvable,
see e.g. Ata and Kumar \cite{AK1} and references therein.

As far as BCP are concerned, the model studied here is a special case of the models considered
in some of the aforementioned papers. In particular, BCP and RBCP
play here important roles,
where the reduction from an $I$-dimensional BCP to a one-dimensional RBCP
is a special case of \cite{Har-Wil}.
Moreover, the HJB equation, that in the present setting is an ordinary differential equation
and will be referred to merely as a {\it Bellman} equation,
is a special case of the partial differential equations treated in \cite{ABW}.
In addition, our specific one-dimensional RBCP, its relation to the Bellman equation,
and its solution go back to Harrison and Taksar \cite{har-tak},
where a singular control problem for a BM is solved.
The solution is given by a reflected BM (RBM), with supporting interval
determined by a free boundary problem associated with the Bellman equation.
This type of free boundary problem first appeared in \cite{har-tak},
and we therefore refer to it as the {\it Harrison-Taksar free boundary problem}.
In our case, the interval is always of the form $[0,{\bf x}^*]$, and we call
${\bf x}^*$ the free boundary point.

On the other hand,
the works \cite{BG2006}, \cite{BG2012} and \cite{AK1}, despite their vast generality,
do not cover the present model as they do not treat admission control and rejection penalties.
Thus, while the BCP is well understood, convergence and
AO issues have not been addressed before. Addressing these issues is the main contribution of this paper.
This is done by proving that the BCP value function constitutes a lower bound on the limit inferior
of QCP costs under any sequence of policies (Theorem \ref{th1}), and then constructing a specific policy
that asymptotically achieves this lower bound (Theorem \ref{th2}).
This AO policy depends explicitly on the system parameters,
except that it also depends on
 the quantity ${\bf x}^*$. Moreover, it uses the well-known {\it $c\mu$ rule} in a novel way, as we now explain.

The structure of the policy alluded to above is simple enough to describe without introducing much notation.
The notation needed is as follows.
For class-$i$ customers, denote holding cost per unit time by $h_i$,
rejection penalty per customer by $r_i$, and reciprocal mean service time
by $\mu_i$. The policy is defined in terms of three elements: The
index $h_i\mu_i$, the index $r_i\mu_i$, and the free boundary point ${\bf x}^*$.
The first index is used for scheduling.
It is precisely the index used for the $c\mu$ priority rule
(a terminology used when $c$, rather than $h$, denotes holding cost per unit time), where
classes are prioritized in the order of $h_i\mu_i$, the highest priority given
to the class $i$ with greatest $h_i\mu_i$.
As observed first by Smith \cite{smith} and Cox and Smith \cite{cox-smi},
the $c\mu$ priority rule is exactly optimal for holding costs.
Many extensions to this result have been
shown (see e.g., \cite{BVW}, \cite{van} and discussions therein).
Our scheduling policy uses the same index to assign priorities, but in a state-dependent fashion, as follows.
At any given time, the lowest priority is assigned to the class $i$
having lowest index $h_i\mu_i$ {\it among classes for which the buffers
are not nearly full}. We give precise meaning to the term `nearly full'.

Let us contrast this with the case of infinite buffers and no rejection.
For this model, an AO policy applying dynamic priorities, in the form of
an extended version of the $c\mu$ rule, was developed by van Mieghem \cite{van}
to address nonlinear delay costs.
When costs are linear, as they are in the present paper,
it is the {\it fixed} priority rule according to $h_i\mu_i$ that is AO.
Suppose now that $I$ is the class that has lowest $h_i\mu_i$ value,
so that class $I$ is assigned lowest priority by this rule.
Then, as is well-known since Whitt \cite{whitt71}, the multiclass G/G/1 queue behaves in such a way that
all classes $i<I$ exhibit vanishing queuelength in the heavy traffic limit.
Consequently, it is not only the aforementioned assignment rule that is AO.
Any priority policy assigning lowest priority to the class $I$
performs equally well, and is therefore AO for such a QCP.
In other words, the only aspect of the index policy which is important for AO in the problem
with infinite buffers and no rejections, is the class assigned the lowest priority.
Thus there is a major difference between the way in which the index is
used in the infinite buffer setting and in this paper.
In the latter case, the full information on the ordering of classes is important.

The admission control is based on the other index,
$r_i\mu_i$, and the free boundary point ${\bf x}^*$.
The significance of this index for admission control in heavy traffic
was first noticed by Plambeck, Kumar and Harrison \cite{PKH}
(see below). Our policy acts as follows.
When the diffusion-scaled workload level exceeds the level ${\bf x}^*$,
all arrivals of one particular class are rejected. This is the class $i$ having the least $r_i\mu_i$ value.
When the workload level is below ${\bf x}^*$, all arrivals are admitted, except rejections
that must take place so as to keep the buffer size constraint valid (namely arrivals that
occur at a time when the corresponding buffer is full). We call these {\it forced rejections}.
A property of the policy that is important for AO is that it maintains, with high probability,
a low number of forced rejections.
As a result, nearly all rejections occur when the workload exceeds ${\bf x}^*$, and only
from one class.
It is to this end that the scheduling policy prioritizes classes with nearly full buffer.

The aforementioned paper \cite{PKH} studies the problem of minimizing rejection
penalties, subject to throughput time constraints, for the multiclass G/G/1 queue in heavy traffic
(see also Ata \cite{ata06} for a closely related formulation).
Each class has a deterministic constraint on the throughput time, and arrivals that are
admitted into the system are assured that, with high probability, their throughput time
constraint will be kept. This property of the policy is referred
to as {\it asymptotic compliance}.
The policy of \cite{PKH} admits all arrivals except
those from the class having lowest $r_i\mu_i$ value,
and only when the workload exceeds a threshold value.
Thus our admission policy resembles that of \cite{PKH},
except that our threshold level is characterized by
the free boundary problem, whereas it is explicit in \cite{PKH}
(their scheduling policy is different than ours).

But the relation of our work to \cite{PKH} is deeper than similarities in the admission policies.
Reiman's snapshot principle \cite{rei-snap},
and the pathwise Little's law, state that, under suitable assumptions,
a deterministic relation holds in the heavy traffic limit between throughput time and queuelength processes.
Accordingly, buffer constraints on queuelength should be asymptotically equivalent
to throughput time constraints.
We follow this rationale in the last section of this paper, where we formulate a QCP that parallels
the QCP addressed in the main body of the paper, where finite buffer constraints are replaced
by throughput time constraints. This may be regarded an extended version of the QCP of \cite{PKH}
that accommodates holding costs.
We do not succeed in fully solving this problem here; our purpose in this part of
the work is mainly to pose the problem and to discuss similarities with the main body
of this paper, leaving the main question open.
We begin by proving a pathwise Little's law
in the form of a conditional result (Proposition \ref{prop3}).
There is no guarantee that
queuelength and throughput times satisfy Little's law under an arbitrary sequence of controls.
We show that $C$-tightness of the processes involved suffices.
Using this result we can show that the
policy we develop for the finite buffer problem satisfies the throughput time constraints and that its
limit performance is dominated by the BCP value (Theorem \ref{th3}).
In order to deduce that it is AO, a lower bound in the same form is also needed.
However, due to the lack of validity of Little's law for general sequences of policies,
we can only show AO in a restricted
class of policies (Proposition \ref{prop4}).
The broader problem, and hence the question of AO remain open (see Conjecture \ref{conj1}).

Under the AO policy, the $I$-dimensional queuelength process converges
to the process solving the BCP. This convergence is a form of a {\it state space collapse} (SSC),
a term referring to a behavior where queuelength process limits are
dictated by workload process limits.
SSC is an important ingredient in the analysis of queueing network
models in heavy traffic. It has been considered in many works, and in particular in
a general setting by Bramson \cite{bramson-ssc} and Williams \cite{williams-ssc}.
The form of the SSC obtained in this paper involves spatial inhomogeneity
due to the dynamic priorities, and is not covered by \cite{bramson-ssc}, \cite{williams-ssc},
or, to the best of our knowledge, any other work on SSC.
A part of the proof of Theorem \ref{th2} is aimed at showing a SSC result.

For a different formulation of a QCP with finite buffers and rejection costs,
see Ghosh and Weerasinghe \cite{GW1}.
For a formulation other than \cite{PKH} that combines asymptotic compliance and asymptotic
optimality see Plambeck \cite{plamb}.
See Ward and Kumar \cite{war-kum}, Rubino and Ata \cite{rub-ata} and Ata and Olsen \cite{ata-ols}
for other treatments of AO in heavy traffic via a Bellman
equation with free boundary, and Dai and Dai \cite{daidai}
for results on heavy traffic for systems with finite buffers
without optimal control aspects.
Finally, see Ghamami and Ward \cite{gha-war}
for asymptotic optimality results based on a Bellman equation
for the BCP, for a model with customer abandonment rather than rejection.

\skp

We will use the following notation. Given $k\in\N$, $\{e^{(i)},i=1,\ldots,k\}$ denote the
standard basis in $\R^k$.
For $x\in\R$, $x^+=\max(x,0)$. For $a,b\in\R^k$, $a=(a_i)_{i=1,\ldots,k}$, $b=(b_i)_{i=1,\ldots,k}$,
we denote $\|a\|=\sum_{i=1}^k|a_i|$ and $a\cdot b=\sum_{i=1}^ka_ib_i$.
For $y:\R_+\to\R^k$ and $T>0$, $\|y\|_T=\sup_{t\in[0,T]}\|y(t)\|$.
The modulus of continuity of $y$ is given by
\begin{equation}\label{104}
\bar w_T(y,\theta)=\sup\{\|y(s)-y(t)\|:s,t\in[0,T],|s-t|\le\theta\}, \qquad \theta,T>0.
\end{equation}
For Polish space $E$,
denote by $\calD_E[0,T]$ the space of RCLL maps from
$[0,T]$ to $E$, equipped with the usual Skorohod topology.
A sequence of stochastic processes with sample paths in this space
is said to be {\it $C$-tight} if it is tight and every
subsequential limit has continuous sample paths w.p.1.
Convergence in distribution of a sequence of
random variables $\{X_n\}$ to $X$ is denoted by $X_n\To X$.
For $a,b\in\R$, $a<b$, $C^2[a,b]$ denotes the set of functions from $[a,b]$ to $\R$
that are twice continuously differentiable on $(a,b)$, for which derivatives of order
$\le 2$ have continuous extensions to $[a,b]$.

The rest of this paper is organized as follows.
In the next section we introduce the queueing model and QCP. Then we formulate the
BCP and the RBCP, and state their solution via the Harrison-Taksar free boundary problem.
We then discuss the interpretation of the solution.
Section \ref{sec3} shows that the BCP value function is a lower bound on
the limit inferior of the sequence of value functions for the QCP.
Section \ref{sec4} constructs a policy for the QCP and proves that it is AO.
Section \ref{sec5} proves pathwise Little's law and
relates the main body of the paper to the throughput time constraints formulation of \cite{PKH}.

\section{Queueing and diffusion models}
\label{sec2}
\beginsec

\subsection{The multiclass G/G/1 model}\label{sec21}

A sequence of systems is considered, indexed by $n\in\N$. Quantities that depend on $n$
have $n$ as superscript in their notation.
The system has a single server and $I\ge1$ buffers, where
each buffer is dedicated to a class of customers.
The capacity of each of the buffers is limited, where the precise formulation of capacity
is presented later.
Customers that arrive at the system may either be accepted or rejected.
Those that are accepted are queued in the corresponding buffers.
Within each class, service is provided in the order of arrival, where
the server only serves the customer at the head of each line.
Processor sharing is allowed, in the sense that the server
is capable of serving up to $I$ customers (of distinct classes) simultaneously.
An {\it allocation vector}, representing the fraction
of effort dedicated to each of the classes, is any member of
\[
\calB:=\Big\{\beta\in\R_+^I: \sum_{i\in\calI}\beta_i\le1\Big\},
\]
where, throughout, $\calI=\{1,2,\ldots,I\}$.

A probability space $(\Om,\calF,\PP)$ is given, on which
all random variables and stochastic processes involved in describing the model
will be defined. Expectation w.r.t.\ $\PP$ is denoted by $\EE$.
Arrivals occur according to independent renewal processes.
Let parameters $\la^n_i >0$, $i\in\calI$, $n\in\N$ be given,
representing the {\it reciprocal mean
inter-arrival times} of class-$i$ customers in the $n$-th system.
Let $\{{\it IA}_i(l) : l\in\N\}_{i\in\calI}$ be independent sequences
of strictly positive i.i.d.\ random variables with mean
$\E[{\it IA}_i(1)]=1$, $i\in\calI$ and squared coefficient of variation
${\rm Var}({\it IA}_i(1))/\E[{\it IA}_i(1)^2]=C^2_{{\it IA}_i}\in(0,\iy)$.
With $\sum_1^0=0$,
the number of arrivals of class-$i$ customers up to time $t$, for the $n$-th system,
is given by
\begin{equation}\label{28}
A^n_i(t)=A_i(\la^n_it), \quad \text{ where }\quad
A_i(t)=\sup\Big\{l\geq 0:\sum_{k=1}^l{\it IA}_i(k)\le t\Big\},\quad t\geq 0.
\end{equation}
The parameters $\la^n_i$ satisfy
\begin{equation}\label{11}
\la^n_i=n\la_i+\sqrt n\hat\la_i+o(\sqrt n),
\end{equation}
where $\la_i>0$ and $\hat\la_i\in\R$ are fixed.

Similarly, let parameters $\mu^n_i>0$, $i\in\calI$, $n\in\N$ be given,
representing {\it reciprocal mean service times} for service to class $i$
in the $n$-th system.
Let independent sequences $\{{\it ST}_i(l) : l\in\N\}_{i\in\calI}$ of strictly
positive i.i.d.\ random
variables (independent of the sequences $\{{\it IA}_i\}$)
be given, with mean $\E[{\it ST}_i(1)]=1$ and squared coefficient of variation
${\rm Var}({\it ST}_i(1))/\E[{\it ST}_i(1)^2]=C^2_{{\it ST}_i}\in(0,\iy)$.
The time required to complete the $l$-th service to
a class-$i$ customer in the $n$-th system is given
by ${\it ST}_i(l)/\la_i^n$ units of time dedicated by the server to
this class. This can otherwise be stated in terms of the
{\it potential service time} processes, given by
\begin{equation}\label{12}
S^n_i(t)=S_i(\mu^n_it), \quad \text{ where }\quad
S_i(t)=\sup\Big\{ l\geq 0:\sum_{k=1}^l{\it ST}_i(k)\leq t\Big\},
\quad t\geq 0.
\end{equation}
$S^n_i(t)$ is the number of class-$i$ jobs completed by the time when the server has dedicated
$t$ units of time to work on jobs of this class.
It is assumed that $\mu^n_i$ satisfy
\begin{equation}\label{14}
\mu^n_i=n\mu_i+\sqrt n\hat\mu_i+o(\sqrt n),
\end{equation}
where $\mu_i>0$ and $\hat\mu_i\in\R$ are fixed.
The first order quantities $\la_i$ and $\mu_i$ are assumed to satisfy the {\it critical
load condition}
\begin{equation}
  \label{15}
  \sum_{i\in\calI}\rho_i=1,\qquad \text{where}\qquad \rho_i=\frac{\la_i}{\mu_i},\,i\in\calI.
\end{equation}

The number of class-$i$ rejections until time $t$
in the $n$-th system is denoted by $Z^n_i(t)$. Since rejections occur only at times
of arrival, we have
\begin{equation}\label{99}
Z^n_i(t)=\int_{[0,t]}z^{n,i}_sdA^n_i(s)
\end{equation}
for some process $z^{n,i}$.

The number of class-$i$ customers present in the $n$-th system at time $t$
is denoted by $X^n_i(t)$.
For simplicity, the initial number of customers, $X^n_i(0)$ is deterministic,
and it is assumed that no partial service has been provided to
any of the jobs present in the system at time zero.
We will call $X^n=(X^n_i)_{i\in\calI}$ the queuelength process.
Let $B^n=(B^n_i)_{i\in\calI}$
be a process taking values in the set $\calB$, representing the
fraction of effort devoted by the server to the various customer classes.
Then
\be\label{2}
T^n_i(t)=\int_0^t B^n_i(s)ds
\ee
gives the time devoted to class-$i$ customers up to time $t$.
The number of service completions of class-$i$ jobs during
the time interval $[0,t]$ can thus be expressed in terms of the potential service process
and the process $T^n_i$ as
\be\label{1}
D^n_i(t)=S^n_i(T^n_i(t)).
\ee
We thus have
\be\label{3}
X^n_i(t)=X^n_i(0)+A^n_i(t)-D^n_i(t)-Z^n_i(t)=X^n_i(0)
+A^n_i(t)-S^n_i(T^n_i(t))-Z^n_i(t),\qquad t\ge0.
\ee
Use the notation $A^n$ for $(A^n_i)_{i\in\calI}$ and similarly for the processes
$S^n$, $Z^n$, $D^n$, $X^n$, $T^n$.
It is assumed that $B^n$ has RCLL sample paths.
By construction, the arrival and potential service processes also have RCLL paths,
and accordingly, so do $D^n$, $Z^n$ and $X^n$.

We define a rescaled version of the processes at diffusion
scale as
\begin{equation}\label{16-}
\hat A^n_i(t)=\frac{A^n_i(t)-\la^n_it}{\sqrt n},\qquad
\hat S^n_i(t)=\frac{S^n_i(t)-\mu^n_it}{\sqrt n},\qquad i\in\calI,
\end{equation}
\[
\hat Z^n(t)=\frac{Z^n(t)}{\sqrt n},\qquad
\hat X^n(t)=\frac{X^n(t)}{\sqrt n}.
\]

We now come to the buffer structure.
A bounded closed convex set with nonempty interior $\calX\subset\R_+^I$ is given,
satisfying $0\in\calX$.
It is assumed that, for every $n$,
the rescaled initial condition $\hat X^n(0)$ lies in $\calX$,
and that the rejection mechanism assures that the buffer constraint is always met, namely:
\begin{equation}
  \label{16}
  \hat X^n(t)\in\calX,\qquad t\ge0,\, a.s.
\end{equation}
For example, the case $\calX=\{y\in\R_+^I:y_i\le b_i,\,i\in\calI\}$
corresponds to a system having a dedicated buffer, of size $b_i\sqrt n$, for each class, $i$.
A single, shared buffer of size $b\sqrt n$ can be modeled by letting
$\calX=\{y\in\R_+^I:\sum_iy_i\le b\}$. In any case, the actual (un-normalized) buffer
size scales like $\sqrt n$.
To meet the constraint \eqref{16}, the control mechanism
must reject some of the arrivals. In particular,
consider a class-$i$ arrival occurring at a time $t$ when
\begin{equation}\label{48}
(X^n(t-)+e^{(i)})/\sqrt n\not\in\calX.
\end{equation}
This arrival has to be rejected so as to keep \eqref{16} valid.
Physically, this situation represents buffers being full, with no available
space to accommodate new arrivals. Such rejections, that occur when \eqref{48} holds,
are often called {\it loss} in the literature.
In our setting, admission/rejection decisions are controlled by the decision maker,
and it is natural to refer to these as part of the rejection control process.
We will refer to them as
{\it forced rejections}, to distinguish them from rejections that occur when
the buffers are not full (i.e., when \eqref{48} does not hold).

The process $U^n=(Z^n,B^n)$ is regarded a control, that is determined based on observations
from the past (and present) events in the system.
The precise definition is as follows.
\begin{definition}\label{def1} {\bf (Admissible control, QCP)}
Fix $n\in\N$ and consider fixed processes $(A^n,S^n)$ given by
\eqref{28} and \eqref{12}. $A^n$ and $S^n$ are called the {\bf primitive processes}.
A process $U^n=(Z^n,B^n)$, taking values in $\R_+^I\times\calB$, having RCLL
sample paths with the processes $Z^n_i$, $i\in\calI$ having nondecreasing sample paths
and given in the form \eqref{99},
is said to be an {\bf admissible control} for the $n$-th system if the following holds.
Let the processes $T^n$, $D^n$, $X^n$ be defined by the primitive and control
processes, $(A^n,S^n)$ and $(Z^n,B^n)$,
via equations \eqref{2}, \eqref{1} and \eqref{3}, respectively.
Then
\begin{itemize}
\item $(Z^n,B^n)$ is adapted to the filtration $\{\calF^n_t\}_{t\ge0}$, where
\(
 \calF^n_t=\sig\{A^n_i(s),D^n_i(s), i\in\calI, s\leq t\};
\)
\item One has a.s., that, for all $i\in\calI$ and $t\ge0$,
\be
X^n_i(t)=0 \quad \text{implies} \quad B^n_i(t)=0.
\label{4}
\ee
\end{itemize}
An admissible control under which the scaled version $\hat X^n$ of $X^n$ satisfies \eqref{16}
is said to {\bf satisfy the buffer constraints}.
\end{definition}
The first bullet above asserts that control decisions are based on the past arrival and
departure events.
The second bullet expresses the fact that jobs from a certain class
can be processed only if there is at least one customer of that class in the system.
We denote the class of all admissible controls $U^n$ by $\tilde\calU^n$,
and the subset of those members of $\tilde\calU^n$ satisfying the buffer constraints,
by $\calU^n$. Except for the last section of this paper, we will refer
to processes in $\calU^n$ as merely admissible controls, for short.
Note that the class $\calU^n$
depends on the processes $A^n$ and $S^n$, but we consider these processes to be fixed.

Fix $\al>0$, $h\in(0,\iy)^I$ and $r\in(0,\iy)^I$. For each $n\in\N$ consider the cost
\be\label{5}
J^n(U^n)=\E\Big[\int_0^\iy e^{-\al t}[h\cdot\hat X^n(t)dt+r\cdot d\hat Z^n(t)]\Big],
\qquad U^n=(Z^n,B^n)\in\tilde\calU^n.
\ee
It will be assumed throughout that, for some $x_0\in\calX$,
\begin{equation}
  \label{20}
  \hat X^n(0)\to x_0,\quad \text{as } n\to\iy.
\end{equation}

The QCP value is given by
\begin{equation}\label{22}
V^n=\inf_{U^n\in\calU^n}J^n(U^n).
\end{equation}
We will be interested in the asymptotic behavior of $V^n$.

Denote by $\theta^n=(\theta^n_i)_{i\in\calI}$ $\theta^n_i=1/\mu^n_i$,
and $\theta=(\theta_i)_{i\in\calI}$, $\theta_i=1/\mu_i$.
The process $\theta^n\cdot X^n$, its
normalized version $\theta^n\cdot\hat X^n$ and the formal limit of the latter,
$\theta\cdot X$, will play an important role in reducing the dimensionality of the
problem. These processes are often referred to as the {\it nominal workload}
(eg., in \cite{PKH}), but we will refer to them simply as {\it workload}.

\subsection{The Brownian control problems}
\label{sec22}

Using \eqref{15}, \eqref{3} and the definition of the rescaled processes, a simple calculation
shows that the following identity holds for $i\in\calI$ and $t\ge0$:
\begin{equation}
  \label{17}
  \hat X^n_i(t)=\hat X^n_i(0)+\hat W^n_i(t)
  +\hat Y^n_i(t)-\hat Z^n_i(t),
\end{equation}
where, denoting $m_i=\hat\la_i-\rho_i\hat\mu_i$,
\begin{equation}
  \label{19}
  m^n_i=\frac{\la^n_i-\rho_i\mu^n_i}{\sqrt n}=m_i+o(1),
\end{equation}
\begin{equation}
  \label{26}
  \hat W^n_i(t)=\hat A^n_i(t)-\hat S^n_i(T^n_i(t))+m^n_it,
\end{equation}
and
\begin{equation}
  \label{18}
  \hat Y^n_i(t)=\frac{\mu^n_i}{\sqrt n}(\rho_it-T^n_i(t)).
\end{equation}
Since $\sum_i\rho_i=1$ and one always has $\sum_iB^n_i(t)\le1$, it follows that
\begin{equation}\label{21}
\theta^n\cdot\hat Y^n \quad \text{is a nonnegative, nondecreasing process}.
\end{equation}

We derive
from \eqref{17}--\eqref{21} and \eqref{5}--\eqref{22} a control problem associated with diffusion
by taking formal limits.
Consider equation \eqref{17}. The scaled initial conditions converge to $x$ by
\eqref{20}. Next,
the centered, rescaled renewal process $\hat A^n_i$ [resp., $\hat S^n_i$]
converges weakly to
a BM starting from zero, with zero mean and diffusion coefficient $\sqrt{\la_i}C_{{\it IA}^i}$
[resp., $\sqrt{\mu_i}C_{{\it ST}^i}$] (see Section 17 of \cite{Bill}).
Heuristically, if the processes involved in \eqref{17} are to give rise
to a limiting BCP
then in particular $\hat Y^n$ are order one as $n\to\iy$.
Thus by \eqref{18} one has that $T^n(t)$ converge to $\rho t$. Thus, taking into
account the time change in the second term of \eqref{26},
$\hat W^n$ is to be replaced a BM starting from zero, with drift vector $m=(m_i)_{i\in\calI}$
and diffusion matrix $\sig=\diag(\sig_i)$, where
\[
\sig_i^2:=\la_iC^2_{{\it IA}^i}+\mu_iC^2_{{\it ST}^i}\rho_i
=\la_i(C^2_{{\it IA}^i}+C^2_{{\it ST}^i}).
\]
Such a process will be called an $(m,\sig)$-BM.
Finally, $\hat Y^n$ gives rise to a process $Y$ for which $\theta\cdot Y$
is nonnegative and nondecreasing, whereas $\hat Z^n$ to a process
having nonnegative, nondecreasing components.

\subsubsection{The BCP}

\begin{definition}\label{def2} {\bf (Admissible control, BCP)}
An {\bf admissible control} for the initial condition $x_0\in\calX$
is a filtered probability space $(\Om',\calF',\{\calF'_t\},\PP')$ for which
there exist an $(m,\sig)$-BM, $W$, and a process $U=(Y,Z)$ taking values in
$(\R_+^I)^2$, with RCLL sample paths, such that the following
conditions hold:
\begin{itemize}
\item
$W$, $Y$ and $Z$ are adapted to $\{\calF'_t\}$;
\item
\begin{equation}\label{29}
\text{For $0\le s<t$, the increment $W(t)-W(s)$ is independent of $\calF'_s$ under $\PP'$;}
\end{equation}
\item
\begin{equation}
  \label{03}
  \text{$\theta\cdot Y$ and $Z_i$, $i=1,\ldots,I$, are nondecreasing;}
\end{equation}
\item With
\begin{equation}
  \label{02}
  X(t)=x_0+W(t)+Y(t)-Z(t),\qquad t\ge0,
\end{equation}
one has
\begin{equation}
  \label{01}
  X(t)\in\calX\qquad \text{for all $t$, $\PP'$-a.s.}
\end{equation}
\end{itemize}
\end{definition}
We write $\calA(x_0)$ for the class of admissible controls for the initial condition $x_0$.
When we write $(Y,Z)\in\calA(x_0)$ it will be understood that these processes carry with them
a filtered probability space and the processes $W$ and $X$.
Moreover, with a slight abuse of notation,
we will write $\EE$ for the expectation corresponding to this probability space.
For $(Y,Z)\in\calA(x_0)$, let
\begin{equation}
  \label{04}
  J(x_0,Y,Z)=\EE\Big[\int_0^\iy e^{-\al t}[h\cdot X(t)dt+r\cdot dZ(t)]\Big].
\end{equation}
The BCP is to find $(Y,Z)$ that minimize $J(Y,Z)$ and achieve
the value
\begin{equation}\label{24}
V(x_0)=\inf_{(Y,Z)\in\calA(x_0)}J(x_0,Y,Z).
\end{equation}

\subsubsection{The RBCP}

The BCP is treated by reduction to a one-dimensional problem.
This is obtained by multiplying equation \eqref{02} and the processes involved in it
by $\theta$. To introduce it,
denote $\bar x_0=\theta\cdot x_0$, $\bar m=\theta\cdot m$ and
\(
\bar\sig^2=\sum\theta_i^2\sig_i^2.
\)
Let
\begin{equation}\label{42}
{\bf x}=\max\{\theta\cdot\xi:\xi\in\calX\}.
\end{equation}
\begin{definition}\label{def3} {\bf (Admissible control, RBCP)}
An {\bf admissible control} for the initial condition $\bar x_0\in[0,\bx]$
is a filtered probability space $(\Om',\calF',\{\calF'_t\},\PP')$ for which
there exist an $(\bar m,\bar\sig)$-BM, $\bar W$, and a process $\bar U=(\bar Y,\bar Z)$
taking values in
$\R_+^2$, with RCLL sample paths, such that the following
conditions hold:
\begin{itemize}
\item
$\bar W$, $\bar Y$ and $\bar Z$ are adapted to $\{\calF'_t\}$;
\item
For $0\le s<t$, the increment $\bar W(t)-\bar W(s)$
is independent of $\calF'_s$ under $\PP'$;
\item
\begin{equation}
  \label{07}
  \text{$\bar Y$ and $\bar Z$ are nondecreasing;}
\end{equation}
\item With
\begin{equation}
  \label{06}
  \bar X(t)=\bar x_0+\bar W(t)+\bar Y(t)-\bar Z(t),\qquad t\ge0,
\end{equation}
one has
\begin{equation}
  \label{05}
  \bar X(t)\in[0,\bx]\qquad \text{for all $t$, $\PP'$-a.s.}
\end{equation}
\end{itemize}
\end{definition}
We write $\bar\calA(\bar x_0)$ for the class of admissible controls for the initial condition
$\bar x_0$. Given $(\bar Y,\bar Z)\in\bar\calA(\bar x_0)$, let
\begin{equation}
  \label{08}
  \bar J(\bar x_0,\bar Y,\bar Z)=
  \EE\Big[\int_0^\iy e^{-\al t}[\bar h(\bar X(t))dt+\bar r d\bar Z(t)]\Big],
\end{equation}
where
\[
\bar h(w)=\min\{h\cdot\xi:\xi\in\calX, \theta\cdot\xi=w\},\qquad w\in[0,{\bf x}],
\]
\[
\bar r=\min\{r\cdot z:z\in\R_+^I,\theta\cdot z=1\}.
\]
Note that $\bar h$ is convex by convexity of the set $\calX$
(in case when $\calX$ is polyhedral, $\bar h$ is also piecewise linear).
Note also that as members of $(0,\iy)^I$, $\theta$ and $h$ cannot be orthogonal,
thus $\bar h(w)>0$ for any $w>0$. Since $\bar h(0)=0$, it follows that $\bar h$ is
strictly increasing. Let
\[
\bar V(\bar x_0)=\inf_{(\bar Y,\bar Z)\in\bar\calA(\bar x_0)}\bar J(\bar x_0,\bar Y,\bar Z).
\]

Toward relating the two problems,
we will need the following additional definitions.
First, the extremal points of the set $\{z\in\R_+^I:\theta\cdot z=1\}$ are precisely
$\theta_i^{-1}e^{(i)}$, namely $\mu_ie^{(i)}$, $i\in\calI$. Hence there exists
(at least one) $i^*$ such that $\zeta^*=\mu_{i^*}e^{(i^*)}$ satisfies
\[
\zeta^*\in\argmin_z\{r\cdot z:z\in\R_+^I,\theta\cdot z=1\}.
\]
Fix such $i^*$ and the corresponding $\zeta^*$.
Note that $i^*$ can alternatively be characterized via
\begin{equation}
  \label{69}
  r_{i^*}\mu_{i^*}=\min_ir_i\mu_i.
\end{equation}
Next, let $\gamma:[0,\bx]\to\calX$ be Borel measurable, satisfying
\begin{equation}\label{36}
\gamma(w)\in\argmin_\xi\{h\cdot\xi:\xi\in\calX,\theta\cdot\xi=w\},\qquad w\in[0,{\bf x}].
\end{equation}
(For the existence of a measurable selection see
Corollary 10.3 in the appendix of \cite{ethkur}).
Note that, by definition,
$\gamma(w)\in\calX$, $\theta\cdot\gamma(w)=w$, and $h\cdot\gamma(w)=\bar h(w)
\le h\cdot\xi$ for every
$\xi\in\calX$ for which $\theta\cdot\xi=w$.
The relation between the problems is as follows.
\begin{proposition}
  \label{prop1}
  Let $x_0\in\calX$ and $\bar x_0=\theta\cdot x_0$.
\\
  i. Given an admissible control $(\Om',\calF',\{\calF'_t\},\PP',W,Y,Z)$ for $x$
  for the (multidimensional) BCP, define $(\bar W,\bar X,\bar Y,\bar Z)$ by
  $(\theta\cdot W,\theta\cdot X,\theta\cdot Y,\theta\cdot Z)$.
  Then $(\bar Y,\bar Z)\in\bar\calA(\bar x_0)$ and $\bar J(\bar x_0,\bar Y,\bar Z)\le
  J(x_0,Y,Z)$.
\\
  ii. Conversely,
  let an admissible control $(\Om',\calF',\{\calF'_t\},\PP',\bar W,\bar Y,\bar Z)$
  for $\bar x_0$
  for the RBCP be given, and assume the probability
  space supports an $(m,\sig)$-BM $W$. Assume $W$ is $\{\calF'_t\}$-adapted
  and satisfies $\theta\cdot W=\bar W$ and \eqref{29}. Construct $(X,Y,Z)$ by
  \begin{equation}\label{71}
  X(t)=\gamma(\bar X(t)),\qquad Z(t)=\zeta^*\bar Z(t),
  \end{equation}
  \begin{equation}\label{72}
  Y(t)=X(t)-x_0-W(t)+Z(t).
  \end{equation}
  Then $(Y,Z)\in\calA(x_0)$, and
  $J(x_0,Y,Z)\le \bar J(\bar x_0,\bar Y,\bar Z)$.
\\
  iii. $V(x_0)=\bar V(\bar x_0)$.
\end{proposition}

\proof
i. We verify that Definition \ref{def3} is satisfied by $(\bar W,\bar X,\bar Y,\bar Z)$.
The first three bullets in that definition are straightforward.
Equation \eqref{06} follows from \eqref{02}, while \eqref{05} from \eqref{01}.

Now, by definition of $\bar h$,
\begin{equation}\label{74}
h\cdot X(t)\ge\bar h(\theta\cdot X(t)),
\end{equation}
and by definition of $\bar r$,
\begin{equation}\label{75}
\int_0^\iy e^{-\al t}r\cdot dZ(t)\ge\int_0^\iy e^{-\al t}\bar r d(\theta\cdot Z(t)).
\end{equation}
Therefore
\[
J(x_0,Y,Z)\ge\bar J(\bar x_0,\bar Y,\bar Z).
\]
\\
ii.
We show that $(Y,Z)\in\calA(x_0)$ by verifying that Definition \ref{def2} is satisfied.
The adaptedness follows by the assumption on $W$ and the construction of $X,Y$ and $Z$.
Property \eqref{29} holds by assumption.
By construction, $(X,Y,Z)$ satisfy \eqref{02}. Property \eqref{01}
holds because, by definition, $\gamma(w)\in\calX$ for all $w\in[0,\bx]$.
$Z_i$ are nonnegative, nondecreasing because so is $\bar Z$, and $\zeta^*_i\ge0$.
Moreover,
\[
\theta\cdot Y(t)=\theta\cdot X(t)-\theta\cdot x_0-\theta\cdot W(t)+\theta\cdot Z(t)
=\bar Y(t).
\]
Hence $\theta\cdot Y$ is nonnegative, nondecreasing.
As a result, $(Y,Z)\in\calA(x_0)$.

Next, note that
\[
h\cdot X(t)=\bar h(\bar X(t)),\qquad r\cdot dZ(t)=\bar rd\bar Z(t).
\]
Therefore
\begin{equation}\label{09}
J(x_0,Y,Z)=\bar J(\bar x_0,\bar Y,\bar Z).
\end{equation}
iii.
The last assertion will follow from the first two once we show that, in (ii), one
can always find $W$ with the stated properties. This is possible by supplementing
the one-dimensional BM $\bar W$ with an $(I-1)$-dimensional BM, independent of $\bar W$,
and augmenting the probability space accordingly.
Specifically, if $\bar W$ is a (one-dimensional) $(\bar m,\bar\sig)$-BM w.r.t.\ a filtration
$\{\bar\calF_t\}_{t\ge0}$ and
$\hat W$ is a standard $(I-1)$-dimensional BM independent of $\bar W$ then
it is not hard to see that an $I\times(I-1)$ matrix $\hat A$ and $I$-dimensional vectors
$\bar A$ and $a$ can be found so that the $I$-dimensional process
$W(t)=\hat A\hat W(t)+\bar A\bar W(t)+at$ is an $(m,\sig)$-BM and one has
$\theta\cdot W(t)=\bar W(t)$, $t\ge0$.
Letting
\[
\calF'_t=\bar\calF_t\vee\sig\{\hat W(s):s\in[0,t]\},\qquad t\ge0,
\]
then gives a filtration with which all conditions of an admissible control for the BCP
are satisfied.
\qed

\subsubsection{The Harrison-Taksar free boundary problem}

The solution to the one-dimensional problem has been studied
by Harrison and Taksar \cite{har-tak} via the Bellman equation.
They showed that the function $\bar V$ is $C^2[0,\bx]$ and solves the equation
\begin{equation}
  \label{25}
  \begin{cases}
  \ds
 \Big[\frac{1}{2}\bar\sig^2f''+\bar mf'-\al f+\bar h\Big]\w f'\w[\bar r-f']=0,
 \qquad \text{in } (0,\bx),
 \\ \\
 f'(0)=0,\qquad
 f'(\bx)=\bar r.
 \end{cases}
\end{equation}
It follows from their work that an optimal control
is one under which the process $\bar X$ is a RBM on a certain subinterval of $[0,\bx]$.
We will consider a RBM as a path transformation of a BM by a {\it Skorohod map},
a map that will later be used in a wider context.
To introduce this map, let $a>0$. The Skohorod map on the interval $[a,b]$,
denoted by $\Gam_{[a,b]}$, is map $D([0,\iy):\R)\to D([0,\iy):\R)^3$.
It is characterized as the solution
map $\psi\to(\ph,\eta_1,\eta_2)$ to the so called {\it Skorohod Problem},
namely the problem of finding, for a given $\psi$, a triplet $(\ph,\eta_1,\eta_2)$, such that
\begin{equation}\label{44}
\ph=\psi+\eta_1-\eta_2,\qquad \ph(t)\in[a,b] \text{ for all } t,
\end{equation}
\begin{equation}\label{45}
\text{$\eta_i$ are nonnegative and nondecreasing, $\eta_i(0-)=0$,
and $\int_{[0,\iy)}1_{(a,b]}(\ph)d\eta_1=\int_{[0,\iy)}1_{[a,b)}(\ph)d\eta_2=0$.}
\end{equation}
By writing $\eta_i(0-)=0$ we adopt the convention that $\eta_i(0)>0$ is regarded a jump
at zero. This convention, in conjunction with $\int_{[0,\iy)}1_{(a,b]}(\ph)d\eta_1=0$
[resp., $\int_{[0,\iy)}1_{[a,b)}(\ph)d\eta_2=0$],
means that if $\psi(0)<a$ [resp., $\psi(0)>b$] then $\ph(0)=a$ [resp., $b$].
If, however,
$\psi(0)\in[a,b]$ then $\ph(0)=\psi(0)$, and $\eta_i$ have no jump at zero.

See \cite{KLRS} for existence and uniqueness of solutions,
and continuity and further properties of the map.
In particular, it is well-known that $\Gam_{[a,b]}$ is continuous in the uniformly-on-compacts
topology.

We now go back to the RBCP. The following is mostly a result of \cite{har-tak}.
\begin{proposition}
  \label{prop2}
  The function $\bar V$ is in $C^2[0,\bx]$ and solves \eqref{25} uniquely among
  all $C^2[0,\bx]$ functions.\footnote{This uniqueness
  question was left open in \cite{har-tak}, at the end of Section 6.}
  Denote $\bx^*=\inf\{y\in[0,\bx]:\bar V'(z)=\bar r \text{ for } z\in[y,\bx]\}$.
  Then $\bx^*\in(0,\bx)$. Fix $\bar x_0\in[0,\bx]$.
  Let $\bar W$ be an $(\bar m,\bar\sig)$-BM and let $\bar X$, $\bar Y$ and $\bar Z$ be
  the corresponding RBM on $[0,\bx^*]$ and boundary terms for $0$ and $\bx^*$, defined as
  \begin{equation}\label{73}
  (\bar X,\bar Y,\bar Z)=\Gam|_{[0,\bx^*]}(\bar x_0+\bar W).
  \end{equation}
  Then $(\bar Y,\bar Z)$ is optimal for $\bar V(\bar x_0)$, i.e.,
  $\bar J(\bar x_0,\bar Y,\bar Z)=\bar V(\bar x_0)$.
\end{proposition}
\begin{remark}\label{rem1}
Note that $\bar X$ has the form $\bar X=\bar x_0+\bar W+\bar Y-\bar Z$.
Moreover, if $\bar x_0>\bx^*$ then $\bar X$ is initially at
$\bx^*$; in particular, $\bar Z(0)=(\bar x_0-\bx^*)^+$.
\end{remark}

\proof
The fact that $\bar V$ is $C^2$ and solves the equation
is proved in \cite{har-tak}, Proposition 6.6 and the discussion that follows.
Uniqueness follows from the uniqueness of solutions in the viscosity sense, for a class
of equations for which the above is a special case \cite{ABW}.
Let us explain how.
It follows from the main result of \cite{ABW} that uniqueness of viscosity
solutions holds for \eqref{25} where the Neumann boundary condition (BC) is replaced by
a state constraint BC (see \cite{ABW} for the definitions of viscosity solutions
and state constraint BC).
As is well-known (and follows directly from the definition),
any $C^2$ function satisfying equation \eqref{25} is also a viscosity solution in the interior
$(0,\bx)$. As for the state constraint BC, it is easy to check (again following directly from
the definition) that any
smooth function satisfying the Neumann BC $f'(0)=0$ and $f'(\bx)=\bar r$,
also satisfies the state constraint BC. This gives the uniqueness.

It is shown in \cite{har-tak} (see the discussion preceding (6.9) therein)
that $f'=\bar r$ on $[y,\bx]$, some $y\in[0,\bx)$. This shows $\bx^*<\bx$.

Next, it is shown in \cite{har-tak} that the control under which $\bar X$ is
a RBM on $[a,\bx^*]$, for some $0\le a<\bx^*$, is optimal.
It remains to show that, in the case considered in this paper, $a=0$. The argument relies on
the fact that $\bar h$ is strictly increasing, as shown in
the discussion following \eqref{08}.

Arguing by contradiction, assume $a>0$. The interval $[a,\bx^*]$ is independent
of the initial condition, and so we are free to choose any $\bar x_0$.
Consider $\bar x_0=0$.
Consider a BM $\bar W$ and the process $\bar X=\bar x_0+\bar W+\bar Y-\bar Z$
that initially has the value $a$ and is given as
a RBM on $[a,\bx^*]$, driven by $\bar W$. By the result of \cite{har-tak}
alluded to above, $(\bar Y,\bar Z)$ is optimal, i.e., $\bar J(0,\bar Y,\bar Z)=\bar V(0)$.
Let $\tau$ be the first hitting time of $\bar X$ at $a+\eps<\bx^*$.
Next, construct on the same probability space another triplet $(\tilde X,\tilde Y,\tilde Z)$,
where $\tilde X$ behaves as a RBM on $[0,\bx^*]$, driven by $\bar W$,
up to the time $\tau$, and starting at
time $\tau$ agrees with $X$ (in particular, it has a jump at time $\tau$).
In other words,
$(\tilde X,\tilde Y,\tilde Z)=(\bar X-a,\bar Y-a,0)$ on $[0,\tau)$,
and $(\tilde X,\tilde Y,\tilde Z)=(\bar X,\bar Y,\bar Z)$ on $[\tau,\iy)$.
Clearly, the cost incurred by $(\tilde X,\tilde Y,\tilde Z)$ on $[\tau,\iy)$,
namely
\[
\int_{[\tau,\iy)}e^{-\al t}[\bar h(\tilde X(t))dt+\bar rd\tilde Z(t)]
\]
is equal to that incurred by $(\bar X,\bar Y,\bar Z)$ on that interval, while, owing
to the strict monotonicity of $\bar h$ and the positivity of $\tau$,
\[
\int_{[0,\tau)}e^{-\al t}\bar h(\tilde X(t))dt<\int_{[0,\tau)}e^{-\al t}\bar h(\bar X(t))dt.
\]
Note that no cost of the form $\bar rd\tilde Z$ is incurred during the time interval
$[0,\tau)$. Taking
expectations shows $\bar J(0,\tilde Y,\tilde Z)<\bar J(0,\bar Y,\bar Z)=\bar V(0)$,
a contradiction. This shows $a=0$.
\qed

As an immediate consequence of the above two results, we obtain an optimal control
for the BCP.
\begin{corollary}
  \label{cor1}
  Let $x_0\in\calX$ and $W$ be an $(m,\sig)$-BM.
  Denote $\bar x_0=\theta\cdot x_0$ and $\bar W=\theta\cdot W$, and let $\bx^*$
  be the free boundary point. Let $(\bar X,\bar Y,\bar Z)$ be defined in terms of
  $\bar W$ as in \eqref{73}, and let $(X,Y,Z)$ be defined in terms of $(W,\bar X,\bar Z)$
  as in \eqref{71}--\eqref{72}. Then $(Y,Z)$ is optimal for $V(x_0)$, namely
  $J(x_0,Y,Z)=V(x_0)$.
\end{corollary}

\subsection{Discussion}\label{sec23}

A brief description of the solution to the BCP is as follows.
The workload process
$\bar X=\theta\cdot X$ is given as a RBM on $[0,\bx^*]$, where the free boundary point $\bx^*$
is dictated by the Bellman equation.
The multidimensional queuelength process $X$ is recovered from $\bar X$ by
\(
X=\gamma(\bar X).
\)
The multidimensional rejection process $Z$ has only one nonzero component,
namely the $i^*$-th component, which increases only when $\bar X\ge\bx^*$.

This structure has an interpretation for the queueing model, that can be used to
identify asymptotically optimal policies. Our main interest will be in the case
of a rectangular domain, namely
\begin{equation}\label{10}
\calX=\{x\in\R^I:0\le x_i\le b_i, i\in\calI\},
\end{equation}
for some fixed $b_i>0$, representing a system where each class
has a dedicated buffer (this will be our assumption in Section \ref{sec4}, although in Section \ref{sec3}
we allow general domains).
In this case, the parameter $\bx$ associated with the RBCP
(defined in \eqref{42}) is given by $\theta\cdot b$.

The BCP solution suggests that, in the queueing model,
rejections should occur only when the scaled workload exceeds the level $\bx^*$,
and only from class $i^*$.
Recall from \eqref{69} that this class is the class for which $r_i\mu_i$ is minimal.
As explained in \cite{PKH},
$i^*$ is the class for which the rejection penalty per unit of work is smallest.

Next, the relation
\begin{equation}\label{90}
\hat X^n=\gamma(\theta\cdot\hat X^n)+o(1)
\end{equation}
between the queuelength and workload processes should hold.
This is a requirement on the scheduling control.
As mentioned in the introduction, when a critically loaded multiclass G/G/1 queue operates
under fixed priority, the queuelength of all classes but one
is asymptotically zero in diffusion scale, the exception being the class
with least priority \cite{whitt71}. This is a simple example of
a scheduling policy that dictates a relation of the form \eqref{90},
where here $\gamma(w)=(0,\ldots,0,w\mu_I)$.
Relation \eqref{90} with a more complicated $\gamma$
appears implicitly when applying the generalized $c\mu$ rule of \cite{van}.
In \cite{atasol} and \cite{atagur} the scheduling policies keep \eqref{90}
where $\gamma$ is a generic minimizing curve.

We can solve for the minimizing curve $\gamma$ in the present setting, where $\calX$ takes the form \eqref{10}.
Equation \eqref{36} can in this case be written as
\[
\gamma(w)\in\argmin_x\{h\cdot x:0\le x_i\le b_i\text{ for all $i$, and } \theta\cdot x=w\},
\quad w\in[0,\bx].
\]
Assume that the classes are labeled in such a way that
\begin{equation}\label{100}
h_1\mu_1\ge h_2\mu_2\ge\cdots\ge h_I\mu_I.
\end{equation}

\begin{figure}[ht]
\begin{center}
\begin{align*}
\includegraphics[width=9em]{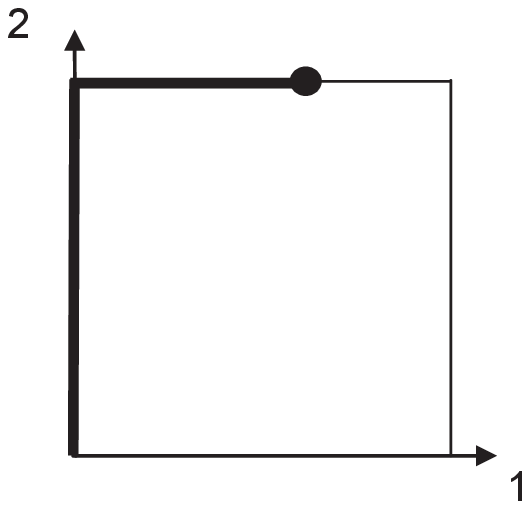}
\hspace{2em}
&
\hspace{2em}
\includegraphics[width=11em]{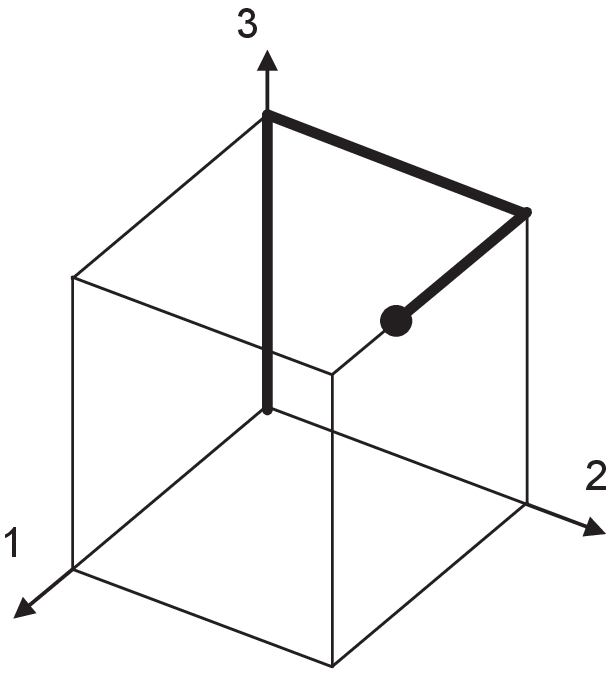}
\end{align*}
\caption{\sl\small
The curve $\gamma$ (thick line)
for the case of dedicated buffers in
dimension $I=2$ and 3. As workload increases, starting from level zero, queuelength $I$
builds up until the corresponding buffer becomes full. Then buffer $I-1$ starts to build up,
and so on,
until the rejection level $(\bx^*,\gamma(\bx^*))$ is reached. Rejections that occur
starting at $\bx^*$ assure that this level is not exceeded.}
\end{center}
\end{figure}

\begin{figure}
\begin{center}
\begin{align*}
\includegraphics[width=9em]{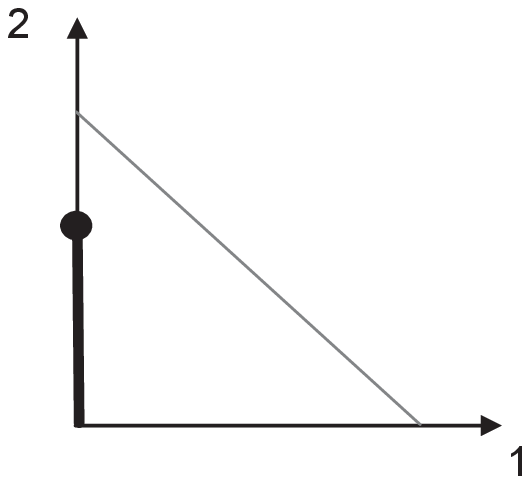}
\hspace{2em}
&
\hspace{2em}
\includegraphics[width=9em]{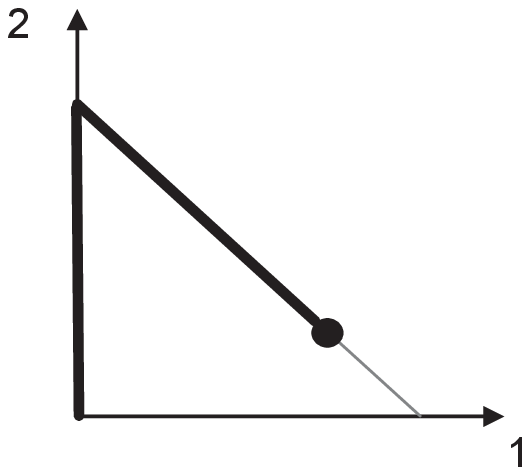}
\end{align*}
\caption{\sl\small
The case of a single buffer shared by 2 classes. Rejection level may be reached before
the buffer with least priority is full (left). For a higher rejection level, the curve continues
along the boundary (right).}
\end{center}
\end{figure}

Given $w\in[0,\bx]$ let $(j,\xi)$ be the unique pair determined by the relation
\begin{equation}\label{37}
w=\sum_{i=j+1}^I\theta_ib_i+\theta_j\xi, \qquad j\in\{1,2,\ldots,I\},\,\xi\in[0,b_j).
\end{equation}
An exception is the special case $w=\bx=\theta\cdot b$, where one lets $j=1$ and $\xi=b_1$.
In other words, denote $\hat b_j=\sum_{i=j+1}^I\theta_ib_i$, $j\in\{0,\ldots,I\}$
and note that $0=\hat b_I<\hat b_{I-1}<\cdots<\hat b_1<\hat b_0=\theta\cdot b=\bx$.
Then $j=j(w)$ is determined by
\[
w\in[\hat b_j,\hat b_{j-1})
\]
and $\xi=\xi(w)=(w-\hat b_j)/\theta_j$. Thus $\gamma$ can be written explicitly as
\begin{equation}\label{38}
\gamma(w)=\sum_{i=j+1}^Ib_ie^{(i)}+\xi e^{(j)}.
\end{equation}

Simple examples are depicted in Figure 1.
While the usual use of the $c\mu$ rule is by assigning fixed priority,
here the index shows up differently.
When buffer $I$ becomes full and workload is increased, a queue in
buffer $I-1$ starts building up, and so on.
A policy that aims at achieving \eqref{90} is developed in Section \ref{sec4}.
Examples for the case of a shared buffer are depicted in Figure 2.
As shown in this figure, the case of two classes with a shared buffer leads to a triangular domain.
In higher dimension one may think of one set of classes
sharing one buffer, another set sharing another
buffer etc., leading to more examples of non-rectangular domains.
General domains are covered in this paper as far as the lower bound is concerned, but we only
address AO controls for the case of rectangular domains.

\begin{example} {\bf (Numerical solution of the BCP)}
  \label{ex1}
In this example we consider a specific three-dimensional BCP and provide its solution explicitly.
The parameters are given in the following table:

\begin{table*}[h]
\begin{center}
\begin{tabular}{c||c|c|c|c|c|c|c}
 $i$ & $b_i$ & $h_i$ & $r_i$
& $\mu_i$ & $\lambda_i$ & $h_i\mu_i$ & $r_i\mu_i$ \\
 \hline
  1 & 15 & 32.9 & 5.0 & 28.0 & 9.33 & 921.2 & 140\\
  2 & 15 & 35.0 & 4.0 & 23.0 & 7.67 &  805 & 92\\
  3 & 10 & 39.0 & 5.5 & 18.0 & 6.0 &  702  & 99\\
\end{tabular}
 \end{center}
 \vspace{-0.3in}
 \label{tab:tb1}
 \end{table*}
 \noi
We further assume that $\hat\la_i=\hat\mu_i=0$ for all $i$ (so that $\bar m=0$), that
$\bar\sig^2=0.1$, and take the discount parameter $\al=10$.
The resulting ordering of the $h_i\mu_i$ index is as in \eqref{100},
namely $h_1\mu_1>h_2\mu_2>h_3\mu_3$. The ordering of $r_i\mu_i$ is such that class 2 is
the most inexpensive as far as rejections are concerned, that is, $i^*=2$.
The Bellman equation takes the form
\begin{equation}
  \label{25n}
  \begin{cases}
  \ds
 [0.05f''-10f+\bar h]\w f'\w[92-f']=0,
 \qquad \text{in } (0,1.74),
 \\ \\
 f'(0)=0,\qquad
 f'(1.74)=92.
 \end{cases}
\end{equation}
The function $\bar h$ defined by
\[
\bar h(w)=\min\Big\{\sum_{i=1}^{3}h_i\xi_i:\xi\in\calX, \sum_{i=1}^{3}\theta_i\xi_i=w\Big\},
\qquad w\in[0,1.74],
\]
where $\calX=[0,15]\times[0,15]\times[0,10]$ and $\theta_i=\mu_i^{-1}$, is
explicitly given by
\[
\bar{h}(w)\approx
\begin{cases}
  \ds
18\cdot 39\, w
 \qquad & 0\leq w\leq0.56,
 \\ \\
\ds
390+23\cdot 35\,(w-0.56)
  \qquad & 0.56<w\leq1.21,
  \\ \\
\ds
915+28\cdot 32.9\,(w-1.21)
    \qquad & 1.21<w\leq1.74.
\end{cases}
\]
A numerical solution of the equation is shown in Figure 3 below.

\begin{figure}
\begin{center}
\begin{align*}
\includegraphics[width=35em, height=16em]{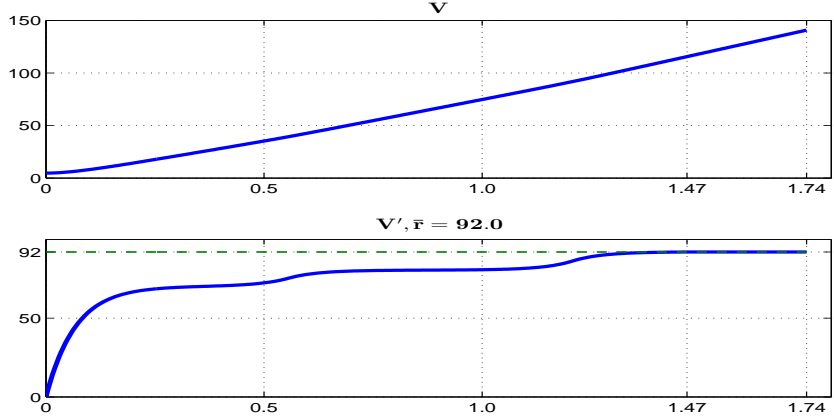}
\end{align*}
\caption{\sl\small
Graphs of $V$ and $V'$. The free boundary point $\bx^*$
is found by seeking the smallest $x$ for which $V'(x)=\bar r$
(Example \ref{ex1}).}
\end{center}
\end{figure}

The free boundary point in this case, found numerically, is the point $\bx^*\approx 1.47$
at which $V'=\bar r$.
The curve $\gamma$ from \eqref{37}--\eqref{38} is given by
\begin{align*}
\gamma(w)\approx
\begin{cases}
  \ds
\frac{w}{0.56}10\, [0, 0, 1]
 \qquad & 0\leq w\leq0.56,
 \\ \\
 \ds
10\, [0, 0, 1] +\frac{w-0.56}{1.21-0.56}15\, [0, 1, 0]
  \qquad & 0.56<w\leq1.21,
  \\ \\
  \ds
10\, [0, 0, 1] +15\, [0, 1, 0] +\frac{w-1.21}{1.74-1.21}15\, [1, 0, 0]
    \qquad & 1.21 < w \leq 1.47.
\end{cases}
\end{align*}
(We have not specified $\gamma$ for values of $w$ beyond the free boundary point $1.47$).
Note that the structure of this curve is of the form depicted in Figure 1 (right).
\end{example}

\begin{example} {\bf (Numerical solution of the QCP)}
  \label{ex2}
Here we present simulation results for the behavior of a two-class M/M/1 queue
operating under the optimal policy.
While for general service time and inter-arrival time distributions finding
the optimal policy is hard, in the case of Poisson arrivals
and exponential service times the problem has the form of a Markov decision process
and one has access to the optimal policy by means of the corresponding Bellman equation
on the discrete 2d grid.
We have solved this equation numerically, computed the optimal policy based on the
solution,
and run a simulation for the behavior of the resulting queuelength process.
Figure 4 depicts histograms for the position
of the two-dimensional queueing process, where gray levels encode the frequency of visits
to each site in the state space (darker gray corresponds to more often visited sites).
The histograms are depicted for an increasing value of the heavy traffic parameter.
The results clearly indicate that the behavior becomes closer and closer
to that given by the limit curve of the form depicted in Figure 1 (left).

While this numerical analysis is related to our results, note carefully that the relation
is indirect: the simulation runs demonstrate the behavior under the optimal policy,
whereas our results address the asymptotics under a sub-optimal (but AO) policy.
The two are related in that both show convergence to the limit behavior identified by
the BCP solution.

Finally, we have also simulated the performance of the sub-optimal policy
that we propose. The graph in Figure 5 shows the ratio between the cost under the
proposed policy and the optimal cost for different values of $n$.

\end{example}

\begin{figure}
\begin{center}
\includegraphics[width=12em, height=12em]{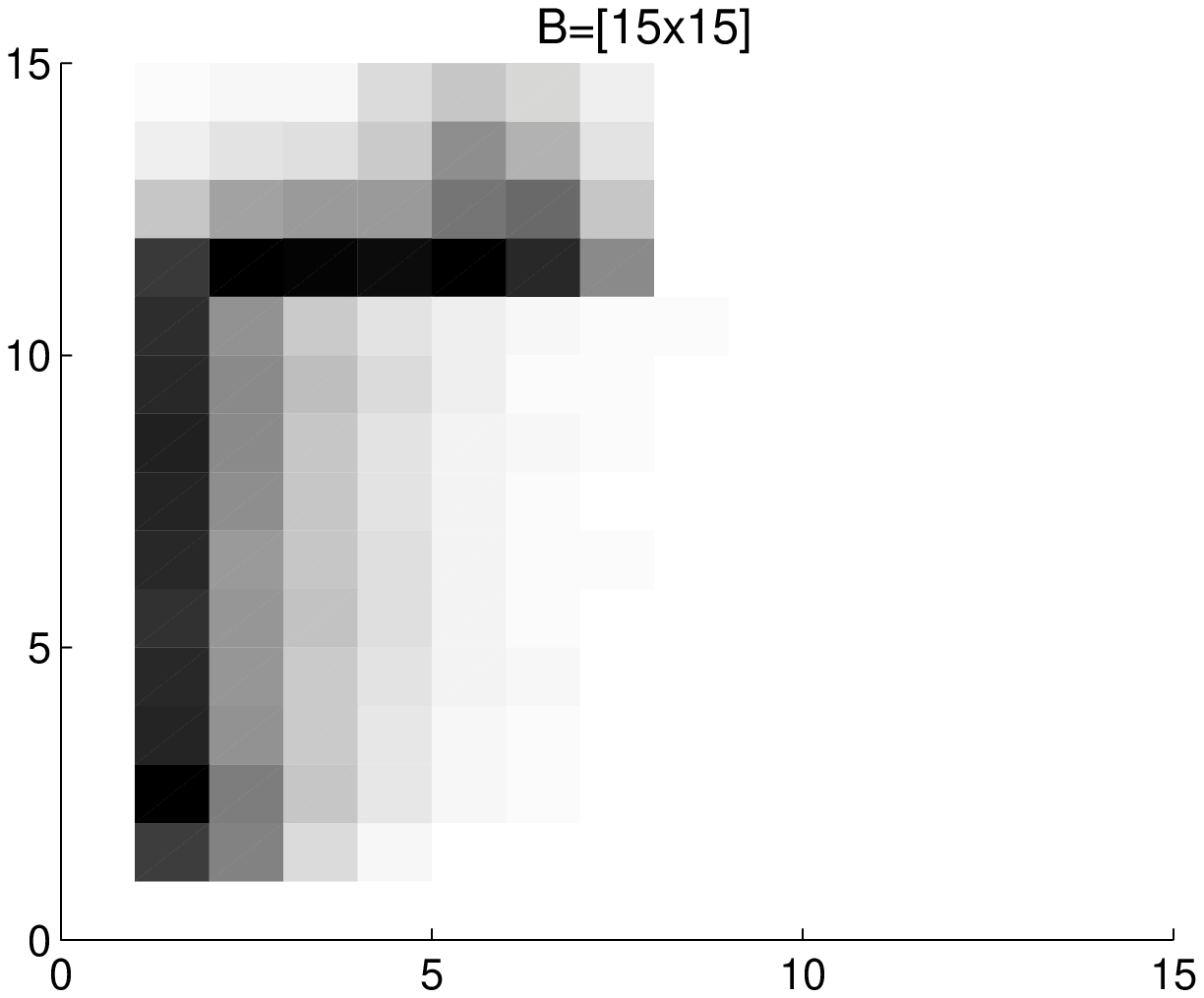}
\hspace{1em}
\includegraphics[width=12em, height=12em]{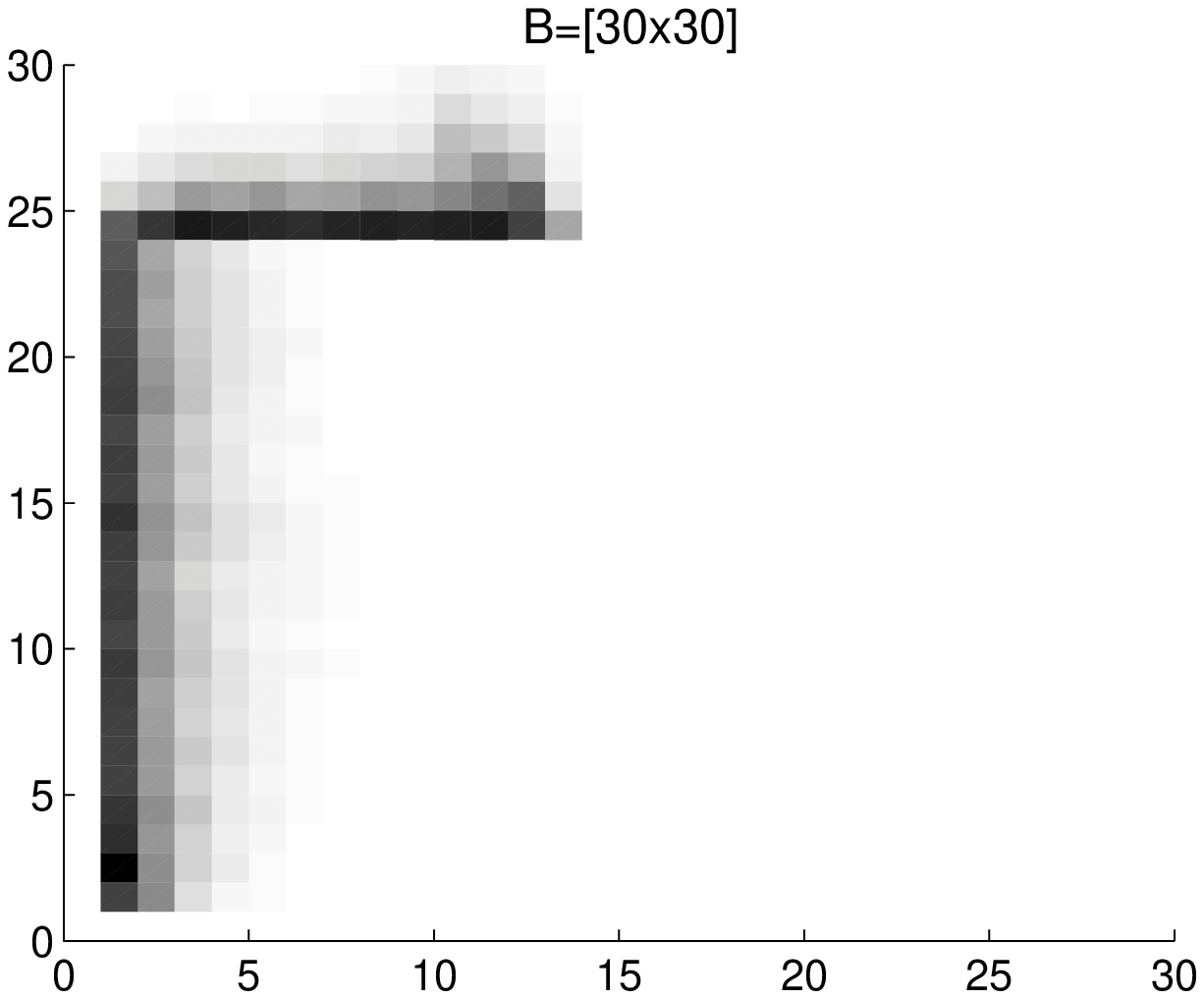}
\hspace{1em}
\includegraphics[width=12em, height=12em]{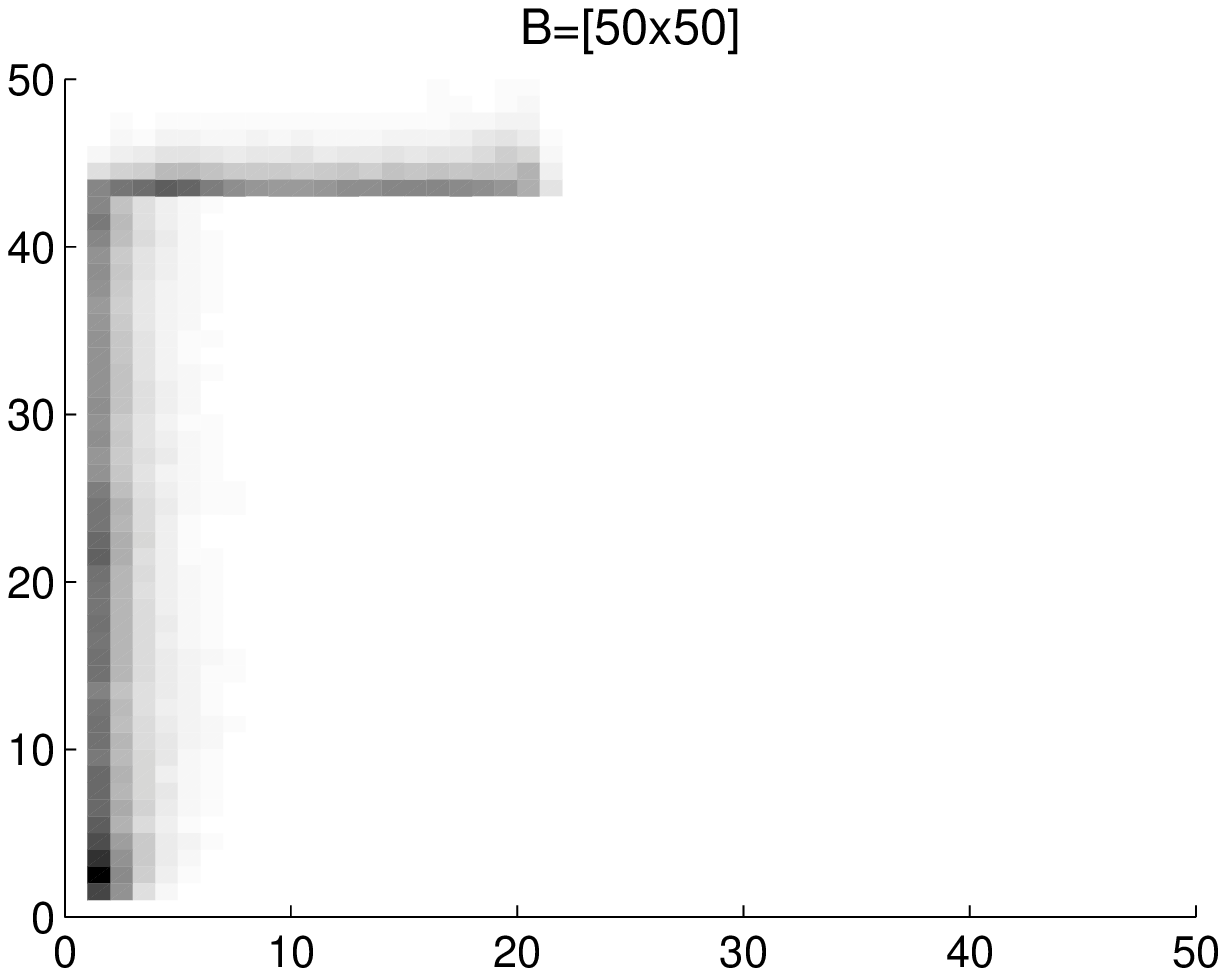}
\caption{\sl\small
Histograms for queuelength process under the optimal policy for an increasing value
of the heavy traffic parameter. The
buffer sizes are $15\times 15$, $30\times30$ and $50\times 50$ (Example \ref{ex2}).}
\end{center}
\end{figure}

\begin{figure}
\begin{center}
\includegraphics[width=22em, height=14em]{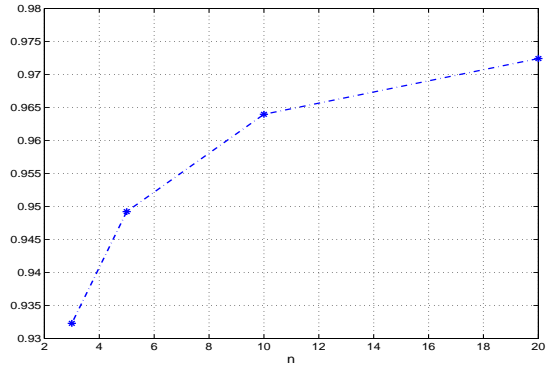}
\caption{\sl\small
Ratio between the (simulated) cost under the proposed policy and the (computed)
optimal cost as a function of $\sqrt{n}$. The graph shows values for $\sqrt{n}=3, 5, 10$
and $20$. The corresponding buffer sizes are given by $5\sqrt{n}$, namely
$15\times 15$, $25\times 25$, $50\times 50$ and $100\times 100$, respectively.
}
\end{center}
\end{figure}

\section{A general lower bound}\label{sec3}
\beginsec

Recall that $V^n$ is defined for the specific initial condition $\hat X^n(0)$,
and that by \eqref{20}, $\hat X^n(0)\to x_0$ as $n\to\iy$. The main result of
this section asserts that the performance of any sequence of policies for the queueing model
is asymptotically bounded below by the BCP value function.
\begin{theorem}
  \label{th1}
  $\uu{V}:=\liminf_{n\to\iy}V^n\ge V(x_0)$.
\end{theorem}
With an eye toward the last section, we will, in fact, prove
a slightly stronger result. Instead of assuming the hard constraint
\eqref{16}, that is a part of the definition of `admissible controls satisfying
the buffer constraint', we will assume throughout this section the following weaker condition.
\begin{equation}\label{110}
\begin{split}
&
\text{For every open set $\tilde\calX\subset\R^I$ with $\calX\subset\tilde\calX$, and every
$T>0$,}
\\
&
\PP(\hat X^n(t)\in\tilde\calX \text{ for all } t\in[0,T])\to1
\text{ as $n\to\iy$.}
\end{split}
\end{equation}

Denote by $\scrI$ the operator
\[
\scrI\ph=\int_0^\cdot\ph(t)dt,
\]
for locally integrable functions $\ph$.

\begin{lemma}\label{lem1}
For $U^n\in\tilde\calU^n$,
\begin{equation}
  \label{27}
  J^n(U^n)=\tilde J^n(U^n):=\E\Big[\int_0^\iy e^{-\al t}
  [\al h\cdot \scrI\hat X^n(t)+\al^2r\cdot\scrI\hat Z^n(t)]dt\Big].
\end{equation}
\end{lemma}
\proof
The identity will follow from integration by parts once we show
that the three terms $e^{-\al t}\scrI\hat X^n(t)$, $e^{-\al t}\hat Z^n(t)$
and $e^{-\al t}\scrI\hat Z^n(t)$
converge to zero a.s.\ as $t\to\iy$.
Note by \eqref{3} that
\[
\hat X^n_i(t)+\hat Z^n_i(t)=n^{-1/2}[X^n_i(t)+Z^n_i(t)]\le n^{-1/2}[X^n_i(0)+A^n_i(t)].
\]
As a renewal process with finite expectation, $A^n_i$ satisfies a law of large numbers
in the sense that $A^n_i(t)/t$ converges a.s.\ as $t\to\iy$.
Thus the three terms alluded to above converge to zero a.s., and the identity follows.
\qed

Before stating the following lemma we introduce some additional notation.
Let $\calA^P(x_0)$ denote the class
of controls for the BCP, defined as in Definition \ref{def2},
except that instead of having RCLL paths, the processes are only assumed to be
progressively measurable. More precisely, an element of $\calA^P(x_0)$ is a filtered
probability space $(\Om',\calF',\{\calF'_t\},\PP')$ with an $(m,\sig)$-BM, $W$,
and a progressively measurable process $(Y,Z)$ taking values in $(R_+^I)^2$, such that
$W$ is adapted, $W(t)-W(s)$ is independent of $\calF'_s$ ($0\le s<t$), and, on an
event having full $\PP'$-measure one has: $\theta\cdot Y$ is a.e.\ equal to a process
with nondecreasing sample paths;
the same holds for each of the processes $Z_i$, $i=1,\ldots,I$;
and, with $X(t)=x_0+W(t)+Y(t)-Z(t)$,
\[
X(t)\in\calX\qquad \text{for a.e.\ $t$.}
\]
Note that $\calA(x_0)\subset\calA^P(x_0)$.

The purpose of introducing this extended class of controls is as follows.
The technique employed in the proof of Theorem \ref{th1} below is based on
tightness of the processes $\scrI\hat X^n$, $\scrI\hat Y^n$ and $\scrI\hat Z^n$
rather than $\hat X^n$, $\hat Y^n$ and $\hat Z^n$. It is established that
the limits of these processes have Lipschitz continuous sample paths, and as a result
they are a.e.\ differentiable. In order to connect these limits to the BCP one needs
to construct from them an admissible control for the latter, but since
the derivatives of Lipschitz functions need not be RCLL, the class of controls $\calA(x_0)$
is too small for this purpose.
Using instead
the class $\calA^P(x_0)$ is possible thanks to a result from \cite{DM-english}
(see below) that shows that
progressively measurable a.e.\ derivatives always exist.

The following lemma shows that working with the extended class of controls does not vary
the value function.

\begin{lemma}\label{lem2}
Let
\[
\tilde J(x_0,Y,Z)=\E\Big[\int_0^\iy e^{-\al t}[\al h\cdot\scrI X(t)+\al^2r\cdot\scrI Z(t)]dt\Big]
\]
and $V^P(x_0)=\inf_{\calA^P(x_0)}\tilde J(x_0,Y,Z)$
(compare with the definitions \eqref{04}, \eqref{24} of $J$ and $V$).
Then $V^P=V$.
\end{lemma}

\proof
Given $x_0$,
consider the specific control that is optimal for $V(x_0)$, namely $(X,Y,Z)$ given
in Proposition \ref{prop1}(ii), where $(\bar X,\bar Y,\bar Z)$ is the RBM on
$[0,\bx^*]$. In particular, $Z(t)=\zeta^*\bar Z(t)$, where $\bar Z$ is one of the boundary
terms of a RBM. It is well known that
$e^{-\al t}(\bar Z(t)+\scrI\bar Z(t))\to0$ a.s., as $t\to\iy$. As a result,
a similar statement holds for $e^{-\al t}(Z_i(t)+\scrI Z_i(t))$, for each $i\in\calI$.
Using integration by parts,
this shows that $V(x_0)=J(x_0,Y,Z)=\tilde J(x_0,Y,Z)\ge V^P(x_0)$.

Next, let $\eps>0$ and consider an $\eps$-optimal control for $V^P(x_0)$, again denoted by
$(X,Y,Z)$. Fix $T>0$.
Construct processes $(\tilde X,\tilde Y,\tilde Z)$ that are identical to $(X,Y,Z)$
on $[0,T)$. As for the time interval $[T,\iy)$,
let $\bar X$ be a RBM on $[0,\bx^*]$ starting from
$\bar X(T)=0$, and let $(\tilde X,\tilde Y,\tilde Z)$ be constructed from
this RBM in the same fashion that $(X,Y,Z)$ are constructed from $\bar X$ in the first
part of the proof. In particular, $\tilde Z$ satisfies
$e^{-\al t}(\tilde Z_i(t)+\scrI\tilde Z_i(t))\to0$,
for each $i\in\calI$ (and it may have a jump at $T$).
By construction, $(\tilde Y,\tilde Z)$ is progressively measurable.
The constructed processes thus form an element of $\calA^P(x_0)$,
and owing to the above tail condition, using integration by parts,
$J(x_0,\tilde Y,\tilde Z)=\tilde J(x_0,\tilde Y,\tilde Z)$. Now,
\begin{align*}
\tilde J(x_0,\tilde Y,\tilde Z)
-\tilde J(x_0,Y,Z)
&=\E\Big[\int_T^\iy
e^{-\al t}[\al h\cdot(\scrI\tilde X(t)-\scrI X(t))
+\al^2r\cdot(\scrI\tilde Z(t)-\scrI Z(t))]dt\Big]\\
&\le c\int_T^\iy te^{-\al t}dt
+\al^2\int_T^\iy e^{-\al t}\int_T^t\E[r\cdot(\tilde Z(s)-\tilde Z(T))]dsdt\\
&\quad
+\E[r\cdot\tilde Z(T)]\al^2\int_T^\iy e^{-\al t}\int_T^tdsdt,
\end{align*}
using the equality $\tilde Z=Z$ on $[0,T)$.
Since on $[T,\iy)$, $\tilde Z-\tilde Z(T)$ is the boundary term of an RBM on a fixed interval,
it is a standard fact that the second term in the above display converges to zero as $T\to\iy$.
As for the last term, since $\tilde J(x_0,\tilde Y,\tilde Z)<\iy$, one has
$\E\int_T^\iy e^{-\al t}r\cdot\scrI\tilde Z(t)dt\to0$ as $T\to\iy$, thus
using monotonicity of $r\cdot\tilde Z$,
\[
e^{-\al (T+2)}\E [r\cdot\tilde Z(T)]\le e^{-\al(T+2)}\E[r\cdot\scrI\tilde Z(T+1)]\le
\E\int_{T+1}^{T+2}e^{-\al t}r\cdot\scrI\tilde Z(t)dt\to0,
\]
as $T\to\iy$.
This shows
\[
J(x_0,\tilde Y,\tilde Z)=\tilde J(x_0,\tilde Y,\tilde Z)\le V^P(x_0)+\eps+a(T),
\]
where $a(T)\to0$ as $T\to\iy$. Taking $T\to\iy$ shows $J(x_0,\tilde Y,\tilde Z)\le V^P(x_0)+\eps$.
Thus to complete the proof, it suffices to show that $J(x_0,\tilde Y,\tilde Z)\ge V(x_0)$.
This is not immediate, because $(\tilde Y,\tilde Z)$ is an element of $\calA^P(x_0)$
whereas $V$ is defined with the smaller class $\calA(x_0)$. We will argue by passing to
the one-dimensional problem. To this end, note that \eqref{74} and \eqref{75} are valid
for the progressively measurable processes, thus
\[
J(x_0,\tilde Y,\tilde Z)\ge\E\Big[\int_0^\iy e^{-\al t}[\bar h(\theta\cdot\tilde X(t))dt+
\bar rd(\theta\cdot\tilde Z(t))]\Big].
\]
Now, the processes $\theta\cdot\tilde Y$ and $\theta\cdot\tilde Z$ are pathwise nondecreasing,
due to the definition of $\calA^P(x_0)$. Hence, if we define $\hat Y(t)=\lim_{s\downarrow t}
\theta\cdot Y(s)$, $\hat Z(t)=\lim_{s\downarrow t}\theta\cdot Z(s)$ and
$\hat X=\theta\cdot x_0+\theta\cdot W+\hat Y-\hat Z$, then $\hat X$, $\hat Y$ and $\hat Z$
are RCLL. Moreover,
they satisfy all assumptions of Definition \ref{def3}, with
$\bar x_0=\theta\cdot x_0$ and $\bar W=\theta\cdot W$.
As a result, they are in $\bar\calA(\bar x_0)$, and so
\[
J(x_0,\tilde Y,\tilde Z)\ge\bar J(\bar x_0,\hat Y,\hat Z)\ge\bar V(\bar x_0).
\]
By Proposition \ref{prop1}, $\bar V(\bar x_0)=V(x_0)$.
We have thus shown that $V(x_0)\le V^P(x_0)+\eps$, and the result follows on taking
$\eps\to0$.
\qed

In the proof below and in the next section we will use
the following characterization of $C$-tightness for processes with sample paths in
$\D_{\R}$ (see Proposition VI.3.26 of \cite{jac-shi}):
{\it
$C$-tightness of $\{X^N\}$, $N\in\N$ is equivalent to}
\begin{align}
&\text{\it The sequence of random variables $\|X^N\|_T$ is tight for every fixed $T<\iy$, and}
\label{105}
\\ \notag \\
&\text{\it For every $T<\iy$,
$\eps>0$ and $\eta>0$ there exist $N_0$ and $\theta>0$ such that} \notag \\
&
N\ge N_0 \text{\it \ implies } P(\bar w_T(X^N,\theta)>\eta)<\eps,
\label{106}
\end{align}
where $\bar w$ is defined in \eqref{104}.

\noi{\bf Proof of Theorem \ref{th1}.}
The structure of the proof is as follows. We invoke
Lemma \ref{lem1} that allows us to work with the cost associated with the integrated
version of the processes. We establish $C$-tightness of the integrated processes;
more precisely, of the sequence $(\hat W^n,\scrI\hat X^n,\scrI\hat Y^n,\scrI\hat Z^n)$.
The rest of the proof is devoted to showing that any subsequential limit of this
sequence gives rise to control within the extended class $\calA^P$,
where the justification to work with the extended class is
provided by Lemma \ref{lem2}.

We thus will rely on Lemma \ref{lem1} and work with $\tilde J^n$.
Using \eqref{22} and Lemma \ref{lem1}, $V^n=\inf_{\calU^n}\tilde J^n(U^n)$.
Fix a subsequence $\{n'\}$ along which $\lim\tilde J^{n'}(U^{n'})=\uu{V}$, and relabel it as $\{n\}$.
Assume, without loss of generality, that $\tilde J^n(U^n)<V(x_0)+1$ for all $n$.
Then $\tilde J^n(U^n)$ is bounded, and so is $J^n(U^n)$, and therefore, for every $T<\iy$,
\[
e^{-\al T}\E[r\cdot\hat Z^n(T)]\le\E\int_0^Tr\cdot d\hat Z^n(t)\le V(x_0)+1.
\]
This shows that $\|\hat Z^n(T)\|$, $n\in\N$, is tight as a sequence of r.v.s,
for each $T$.

Recall that $\hat A^n$ and $\hat S^n$ converge u.o.c.\ to BMs, and note by \eqref{2} that
$T^n_i(t)\le t$ for every $t$. Using this and equations
\eqref{19} and \eqref{26} shows that the sequence of processes $\hat W^n$ is $C$-tight.

Given $T$, using the monotonicity of $\hat Z^n_i(\cdot)$, the Lipschitz constant of
$\scrI\hat Z^n_i|_{[0,T]}$ is bounded by $\|\hat Z^n(T)\|$.
Thus, using the characterization \eqref{105}--\eqref{106},
the tightness of $\hat Z^n(T)$ for each $T$
implies that $\scrI\hat Z^n$ is a $C$-tight sequence of processes.
The condition \eqref{110} implies that, for every $T$, $\|\hat X^n\|_T$, $n\in\N$,
is a tight sequence of random variables.
As a result, by \eqref{105}--\eqref{106},
the sequence $\scrI\hat X^n$ is also $C$-tight. Next, by \eqref{17},
\begin{equation}
  \label{30}
  \|\hat Y^n(t)\|\le \|\hat X^n(t)\|+\|\hat W^n(t)\|+\|\hat Z^n(t)\|.
\end{equation}
It follows from this discussion that, for each $T$,
\[
L_n(T):=\|\hat X^n\|_T\vee\|\hat Y^n\|_T\vee\|\hat Z^n(T)\|
\]
is a tight sequence of r.v.s, and that $(\scrI\hat X^n,\scrI\hat Y^n,\scrI\hat Z^n)$
is $C$-tight, with bound $L_n(T)$ on the Lipschitz constant over the interval
$[0,T]$. Since $L_n(T)$
are tight for each $T$, any weak limit point of the $C$-tight sequence is a process
having locally Lipschitz paths a.s.

Next, since for each $T$, the sequence $\|\hat Y^n\|_T$ is
tight, we have by \eqref{18} and the fact $\mu^n_i/\sqrt n\to\iy$, that
$T^n_i$ converge u.o.c.\ to $\bar T_i$ where $\bar T_i(t)=\rho_it$.
By \eqref{26}, using a lemma regarding random change of
time \cite{Bill}, p.\ 151, it follows that $\hat W^n\To W$, where we recall
that $W$ is an $(m,\sig)$-BM.

By tightness of
$(\hat W^n,\scrI\hat X^n,\scrI\hat Y^n,\scrI\hat Z^n)$,
there exists a convergent subsequence. Denote its limit by
$(W,IX,IY,IZ)$. Note that the last three terms have Lipschitz sample paths.
By an argument as in section IV.17 of \cite{DM-english}, they possess a.e.\ derivatives that are progressively
measurable w.r.t.\ the filtration $\calF'_t=\sig\{W(s),IX(s),IY(s),IZ(s):s\le t\}$.
For concreteness, let $\pl^-f$, for a Lipschitz $f:[0,\iy)\to\R$,
be defined by $\pl^-f(0)=0$ and
$\pl^-f(t)=\liminf_{s\uparrow t}(t-s)^{-1}(f(t)-f(s))$, $t>0$.
Define pathwise $X=(X_i)$, $Y=(Y_i)$
and $Z=(Z_i)$ as $X_i=\pl^-IX_i$, $Y_i=\pl^-IY_i$ and $Z_i=\pl^-IZ_i$.
Then $(X,Y,Z)$ are progressively measurable, and $\scrI X=IX$.
We will show below that these processes along with the filtration $\{\calF'_t\}$
form an element of the class $\calA^P(x_0)$.
Consequently, using Lemma \ref{lem1} and Fatou's lemma for the subsequence under consideration,
\[
\uu{V}=\liminf\tilde J^n(U^n)\ge\tilde J(x_0,Y,Z)
\ge\inf_{\calA^P(x_0)}\tilde J(x_0,\cdot,\cdot)=V^P(x_0)
=V(x_0),
\]
where in the last equality we used Lemma \ref{lem2}.

It thus remains to show that the progressively measurable processes
we have constructed form an element of the class $\calA^P(x_0)$.
To show \eqref{29}, we borrow a few lines from the proof of Lemma 6 of \cite{AMR}.
Fix $0\le s\le t<t+u$. Let
$\al^n=(\hat W^n(s),\scrI \hat X^n(s),\scrI \hat Y^n(s),\scrI \hat Z^n(s))$
and $\al=(W(s),IX(s),IY(s),IZ(s))$.
For $i\in\calI$ let $t^n_i$ [resp., $\tau^n_i$] denote the
renewal epoch of $A^n_i$ [resp., $S^n_i$] following $t$ [resp., $T^n_i(t)$]. That is,
\[
t^n_i=\inf\{t'\ge t:A^n_i(t')>A^n_i(t)\},\qquad
\tau^n_i=\inf\{t'\ge T^n_i(t):S^n_i(t')>S^n_i(T^n_i(t))\}.
\]
Let $\beta^n=(\beta^n_i)_{i\in\calI}$ be defined by
\[
\beta^n_i=(A^n_i(t^n_i+u)-A^n_i(t^n_i),\,S^n_i(\tau^n_i+\rho_iu)-S^n_i(\tau^n_i)).
\]
Then $\al^n$ and $\beta^n$ are mutually independent.
As a result, $\al^n$ and $\gamma^n=(\gamma^n_i)_{i\in\calI}$ are mutually independent, where
\[
\gamma^n_i=\hat A^n_i(t^n_i+u)-\hat S^n_i(\tau^n_i+u)-\hat A^n_i(t^n_i)
+\hat S^n_i(\tau^n_i)+m^n_iu.
\]
Recall the definition \eqref{26} of $\hat W^n$.
We have $t^n_i\To t$, and $T^n(t)\To\bar T(t)$ by which $\tau^n_i\To\rho_it$.
As a result, $\hat W^n_i(t+u)-\hat W^n_i(t)-\gamma^n\To0$. This shows that
$\al$ and $W(t+u)-W(t)$ are mutually independent.
Since $u>0$ and $s\le t$ are arbitrary, an application of
Theorem 1.4.2 of \cite{durrett-book}
shows that all increments $W(t+u)-W(t)$ and $\calF'_s$ are independent.

Let $\calX_\del=\{x\in\R^I:\dist(x,\calX)<\del\}$, $\del>0$.
Condition \eqref{110} implies that for every $s<t$ and $\del>0$,
$(t-s)^{-1}(\scrI\hat X^n(t)-\scrI\hat X^n(s))\in\calX_\del$ occurs with probability tending to $1$
as $n\to\iy$.
As a result, $(t-s)^{-1}(IX(t)-IX(s))\in\oo{\calX_\del}$ a.s. Since $\calX$ is closed and convex,
the intersection of $\oo{\calX_\del}$ over $\del>0$ gives $\calX$,
so $(t-s)^{-1}(IX(t)-IX(s))\in\calX$ a.s.
Thus $X(t)\in\calX$ for a.e.\ $t$, a.s.
Now, each $\scrI\hat Z^n_i$ is nonnegative, nondecreasing and convex, hence so is
$IZ_i$. Therefore $Z_i$ is nonnegative and nondecreasing. As for
$\theta\cdot Y$, note that it is a.e. equal to the pathwise left-derivative of the
process $\theta\cdot IY$, which, for reasons as above, has convex sample paths a.s. Hence $\theta\cdot Y$
is a.e.\ equal to a nondecreasing process.
This shows that $(X,Y,Z)\in\calA^P(x_0)$ and completes the proof.
\qed

\section{A nearly optimal policy in the case of a rectangle}\label{sec4}
\beginsec

In this section we consider the case of a rectangular domain, where
each customer class has a dedicated buffer.
We have introduced in Section \ref{sec23} some notation for this case,
and identified the curve $\gamma$. In particular, the domain $\calX$ is given by
\eqref{10},
where $b_i>0$ are fixed constants, and
the parameter $\bx$ is given by $\theta\cdot b$.
The classes are labeled so that
\[
h_1\mu_1\ge h_2\mu_2\ge\cdots\ge h_I\mu_I,
\]
and, given $w\in[0,\bx]$, $(j,\xi)=(j,\xi)(w)$ are determined by
\[
w\in[\hat b_j,\hat b_{j-1})
\]
and $\xi=\xi(w)=(w-\hat b_j)/\theta_j$,
where $\hat b_j=\sum_{i=j+1}^I\theta_ib_i$, $j\in\{0,\ldots,I\}$
and one has $0=\hat b_I<\hat b_{I-1}<\cdots<\hat b_1<\hat b_0=\theta\cdot b=\bx$.
With this notation, $\gamma$ is given (as in \eqref{38}) by
\[
\gamma(w)=\sum_{i=j+1}^Ib_ie^{(i)}+\xi e^{(j)}.
\]

The difficulty in treating the queueing model according to the BCP solution, as described
in terms of $\gamma$, is that this curve lies along the boundary of $\calX$,
in particular,
along the part $\pl^+\calX:=\{x\in\calX:x_i=b_i \text{ for some } i\}$
of the boundary $\pl\calX$. This part corresponds to states at which
some of the buffers are full. This sets up contradictory goals of keeping some of the
buffers (nearly) full and at the same time avoiding any rejections except when
the workload process reaches the level $\bx^*$.
The policy we propose is based on an approximation of $\gamma$ by another curve
that is bounded away from the buffer limit boundary.

Let $\eps\in(0,\min_ib_i)$ be given.
Let $a_i=b_i-\eps$, $i\in\calI$, and
$a^*:=\bx^*\w(\theta\cdot a)<\bx=\theta\cdot b$. Note that if $\eps$ is small
then $a^*=\bx^*$ (unless $\bx^*=\bx$).
We define an approximation $\gamma^a:[0,\bx]\to\calX$ of $\gamma$
by first defining it on $[0,\theta\cdot a]$ as the function obtained upon replacing
the parameters $(b_i)$ by $(a_i)$ in \eqref{37} and \eqref{38}.
That is, for $w\in[0,\theta\cdot a)$, the variables $j=j(w)$ and $\xi=\xi(w)$ are determined via
\begin{equation}\label{61}
w=\sum_{i=j+1}^I\theta_ia_i+\theta_j\xi, \qquad j\in\{1,2,\ldots,I\},\,\xi\in[0,a_j),
\end{equation}
and
\begin{equation}\label{62}
\gamma^a(w)=\sum_{i=j+1}^Ia_ie^{(i)}+\xi e^{(j)}.
\end{equation}
Given $w\in[0,\theta\cdot a)$, we will sometimes refer to the unique pair
$(j,\xi)$ alluded to above as the {\it representation $(j,\xi)$ of $w$ via \eqref{61}}.
Next, on $[\theta\cdot a,\theta\cdot b]$ we only need the function $\gamma^a$
to be continuous and satisfy the relation
$\theta\cdot\gamma^a(w)=w$. For concreteness we may define
it as the linear interpolation between the points $(\theta\cdot a,a)$
and $(\theta\cdot b,b)$:
\[
\gamma^a(w)=a+\frac{w-\theta\cdot a}{\theta\cdot b-\theta\cdot a}(b-a),
\qquad w\in[\theta\cdot a,\theta\cdot b].
\]
We also define $\hat a_j=\sum_{i=j+1}^I\theta_ia_i$, $j\in\{0,1,\ldots,I\}$,
similarly to $\hat b_j$.

The definition of the policy is provided by specifying $(Z^n(t),B^n(t))$
as a function of $X^n(t)$.

{\it Rejection policy:}
As under any policy, in order to keep the
buffer size constraint \eqref{16}, all forced rejections take place. That is,
if a class-$i$ arrival occurs at a time $t$ when $\hat X^n_i(t-)+n^{-1/2}>b_i$,
then it is rejected.
Apart from that, no rejections occur from any class except class $i^*$,
and no rejections occur (from any class) when $\theta\cdot \hat X^n<a^*$.
When $\theta\cdot \hat X^n\ge a^*$, all class-$i^*$ arrivals are rejected.

\begin{figure}
\begin{center}
\begin{align*}
\includegraphics[width=14em]{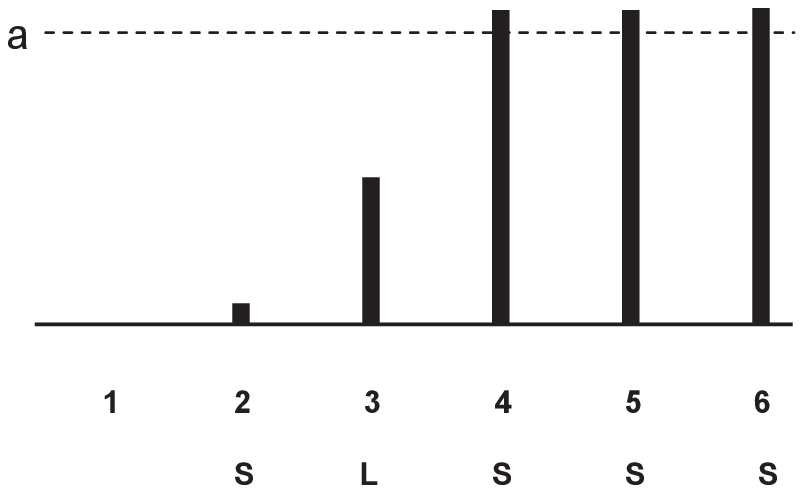}
\hspace{2em}
&
\hspace{2em}
\includegraphics[width=14em]{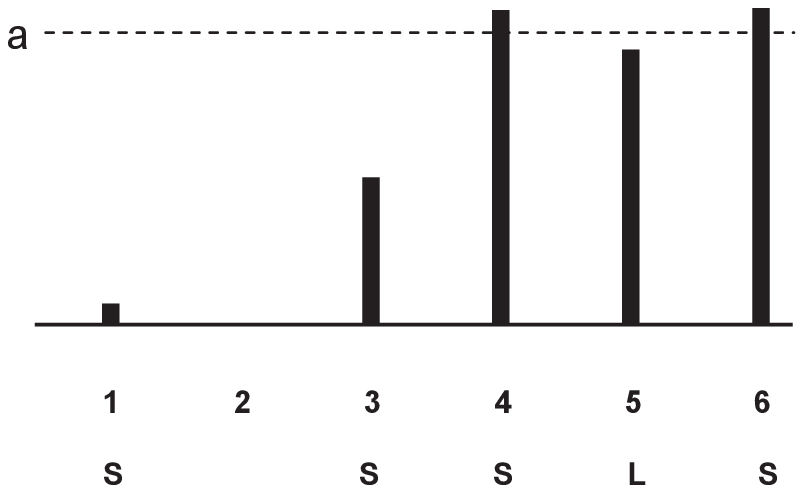}
\end{align*}
\caption{\sl\small
An example with imaginary buffer sizes $a_i=1$
and reciprocal service rates $\theta_i=1$, $i=1,\ldots,6$. The figures depict possible states
$\hat X^n(t)=x$
at a time when the normalized workload $\theta\cdot x=x_1+\cdots+x_6$ is around $3.5$.
The target population distribution is then $\gamma(3.5)=(0,0,0.5,1,1,1)$.
The class of low priority (L) is the maximal $i$ with $x_i<a_i$.
The classes served (S) are the high priority classes having positive population.
Thus, for all $i$, $i$ is being served provided that
$x_i$ exceeds the target population $\gamma_i(3.5)$.
}
\end{center}
\end{figure}

{\it Service policy:}
For each $x\in\calX$ define the class of low priority
\[
\calL(x)=\max\{i:x_i<a_i\},
\]
provided $x_i<a_i$ for some $i$, and set $\calL(x)=I$ otherwise.
The complement set is the set of high priority classes:
\[
\calH(x):=\calI\setminus\{\calL(x)\}.
\]
When there is at least one class among $\calH(x)$ having at least one customer in the system,
$\calL(x)$ receives no service, and
all classes within $\calH(x)$, having at least one customer,
receive service at a fraction proportional to their
traffic intensities. More formally, denote $\calH^+(x)=\{i\in \calH(x):x_i>0\}$,
and define $\rho'(x)\in\R^I$ as
\begin{equation}\label{58}
\rho'_i(x)=\begin{cases}
0, & \text{if } x=0,\\
\ds\frac{\rho_i1_{\{i\in\calH^+(x)\}}}{\sum_{k\in\calH^+(x)}\rho_k}, &
\text{if } \calH^+(x)\ne\emptyset,\\
e^{(I)},& \text{if $x_i=0$ for all $i<I$ and $x_I>0$.}
\end{cases}
\end{equation}
(Note that $\calH^+(x)=\emptyset$ can only happen if $x_i=0$ for all $i<I$,
which is covered by the first and last cases in the above display).
Then for each $t$,
\begin{equation}
  \label{32}
  B^n(t)=\rho'(\hat X^n(t)).
\end{equation}
Note that when $\calH^+(x)\ne\emptyset$,
\begin{equation}
  \label{31}
  \rho'_i(x)>\rho_i\quad \text{for all } i\in \calH^+(x).
\end{equation}
That is, all prioritized classes receive
a fraction of effort strictly greater than the respective traffic intensity.
Also note that $\sum_iB^n_i=1$ whenever $\hat X^n$ is nonzero. This is therefore
a work conserving policy. See Figure 3 for an example of how the class with low priority
and the served classes are determined.
\begin{remark}
The only properties from the structure \eqref{58}--\eqref{32}
that are actually used in the proof are \eqref{31} and $\sum_iB^n_i=1$ if $\hat X^n\ne0$;
in other words we could have
allowed other choices of $\rho'$ as long as \eqref{31} and the work conserving
property hold.
\end{remark}
\begin{remark}
  \label{rem3}
  Although we have assumed $\eps>0$, the policy is well-defined even for $\eps=0$,
  in which case $a=b$, $a^*=\bx^*$ and $\gamma^a=\gamma$. This policy
  is not used here but it is used in the next section.
\end{remark}

Arguing by induction on the times when the driving processes $A^n$ and $S^n$
jump, it is clear that there exists a unique solution to the set of equations
\eqref{99}--\eqref{3}, \eqref{32} along with the verbal description of the rejection
mechanism. Thus the policy is well-defined.

\begin{theorem}
  \label{th2}
  For each $\eps>0$ and $n$, denote the policy constructed above by $U^n(\eps)$. Then
  $\limsup_{n\to\iy}J^n(U^n(\eps))\le V(x_0)+\al(\eps)$, where $\al(\eps)\to0$ as $\eps\to0$.
\end{theorem}
\begin{remark}
i. By a usual diagonalization argument one can extract from $(U^n(\eps),\eps)$
a sequence $U^n$ that is asymptotically optimal, i.e.,
$\limsup_{n\to\iy}J^n(U^n)=V(x_0)$.
\\
ii. The combination of Theorems \ref{th1} and \ref{th2} gives $\uu{V}=V(x_0)$.
\end{remark}

\proof
We fix $\eps$ and write $U^n=(Z^n,B^n)$ for $U^n(\eps)$.
We denote by $\tau^n$ the time of the first forced rejection. A crucial point
about the proof idea is that most of the analysis is performed on the processes
up to the first forced rejection. It is established that the target state
is asymptotically achieved
by the proposed policy, in the sense of weak convergence as $n\to\iy$.
This is done in two steps: First, the workload process
$\theta^n\cdot\hat X^n$ is shown to converge to a RBM,
and then it is shown that $\hat X^n$ lies close to the minimizing curve at all times.
Once these elements are established,
it follows that in any finite time, $\tau^n$ is not reached, and as a result
one has that (i) only rejections from class $i^*$ occur, and only when
$\theta^n\cdot\hat X^n\approx a^*$; (ii) the running cost is minimized locally.
These elements are then combined with some integrability conditions at the last step
of the proof.

We begin with the case where the system starts with initial condition
close to the minimizing curve. More precisely,
\begin{equation}\label{49}
\hat X^n(0)-\gamma^a(\theta\cdot\hat X^n(0))\to0 \text{ as $n\to\iy$,
and }\theta^n\cdot\hat X^n(0)\in[0,a^*] \text{ for all $n$ large.}
\end{equation}
At the last step of the proof we relax this assumption.

{\it Step 1. $C$-tightness for the workload and related processes.}
We multiply equation \eqref{17} by the vector $\theta^n=(1/\mu^n_i)_{i\in\calI}$ and
denote
\begin{equation}
  \label{33}
  W^{\#,n}=\theta^n\cdot\hat W^n,\quad
  X^{\#,n}=\theta^n\cdot\hat X^n,\quad
  Y^{\#,n}=\theta^n\cdot\hat Y^n,\quad
  Z^{\#,n}=\theta^n\cdot\hat Z^n.
\end{equation}
We have
\begin{equation}
  \label{34}
  X^{\#,n}=X^{\#,n}(0)+W^{\#,n}+Y^{\#,n}-Z^{\#,n}.
\end{equation}
Let $W^{\circ,n}:=W^{\#,n}(\cdot\w\tau^n)$ denote the process $W^{\#,n}$ when
stopped at the time $\tau^n$. Define similarly $X^{\circ,n}$, $Y^{\circ,n}$
and $Z^{\circ,n}$. Our goal in the step is to show that the sequence
$(W^{\circ,n},X^{\circ,n},Y^{\circ,n},Z^{\circ,n})$ is
$C$-tight, and that any subsequential limit $(\tilde W,\tilde X,\tilde Y,\tilde Z)$
satisfies a.s.,
\begin{equation}\label{43}
(\tilde X,\tilde Y,\tilde Z)=\Gam_{[0,a^*]}[\bar x_0+\tilde W].
\end{equation}

To this end, note first that the argument for $C$-tightness of the processes $\hat W^n$,
given in the proof of the lower bound,
is valid here. As a result, $W^{\#,n}$ are $C$-tight. Hence so are $W^{\circ,n}$.

By construction (see \eqref{32}), the policy is work conserving, namely
$\sum_iB_i^n(t)=1$ whenever $\hat X^n(t)$ is nonzero.
By the relations \eqref{2} and \eqref{18}, it follows that
the nondecreasing process $\hat Y^{\#,n}$ does not increase when
$X^{\#,n}>0$. A similar property then holds for the stopped processes,
and this can be expressed as
\begin{equation}\label{47}
\int 1_{\{X^{\circ,n}(t)>0\}}dY^{\circ,n}(t)=0.
\end{equation}
Fix $T>0$.
We show next that, as $n\to\iy$,
\begin{equation}\label{46}
\Big(\sup_{t\in[0,T]}X^{\circ,n}(t)-a^*\Big)^+\To0.
\end{equation}
For $\eps'>0$ consider the event $\Om^n_1:=\{\sup_{t\in[0,T]}X^{\circ,n}(t)>a^*+\eps'\}$.
On this event there exist random times $0\le\tau^n_1<\tau^n_2\le\tau^n$ such that
$X^{\#,n}(\tau^n_1)\le a^*+\eps'/2$, $X^{\#,n}(\tau^n_2)\ge a^*+\eps'$ and
$X^{\#,n}(t)>a^*$ for all $t\in[\tau^n_1,\tau^n_2]$.
Thus by \eqref{34} and the fact that $Y^{\#,n}$ does not increase on an interval
where the system is not empty,
denoting here and in the sequel $A[s,t]=A(t)-A(s)$ for any process $A$,
\begin{align*}
(a^*+\eps')-(a^*+\eps'/2) &\le X^{\#,n}[\tau^n_1,\tau^n_2]\\
&=W^{\#,n}[\tau^n_1,\tau^n_2]-Z^{\#,n}[\tau^n_1,\tau^n_2]\\
&=W^{\#,n}[\tau^n_1,\tau^n_2]
-\frac{A^n_{i^*}[\tau^n_1,\tau^n_2]}{\sqrt n}
\end{align*}
where we used the fact that the policy rejects all class-$i^*$ jobs
when $X^{\#,n}>a^*$.
Fix a sequence $r_n>0$, $r_n\to0$, such that $\sqrt n r_n\to\iy$.
In case $\tau^n_2-\tau^n_1<r_n$, the above implies
\[
\eps'/2\le W^{\#,n}[\tau^n_1,\tau^n_2]\le\bar w_T(W^{\#,n};r_n).
\]
In case $\tau^n_2-\tau^n_1\ge r_n$,
\[
2\|W^{\#,n}\|_T\ge\frac{A^n_{i^*}[\tau^n_1,\tau^n_2]}{\sqrt n}
=\hat A^n_{i^*}[\tau^n_1,\tau^n_2]
+\frac{\la^n_{i^*}}{\sqrt n}(\tau^n_2-\tau^n_1)
\ge -2\|\hat A^n_{i^*}\|_T+c\sqrt nr_n,
\]
for some positive constant $c$. Combining the two cases, the $C$-tightness of
$W^{\#,n}$ and the tightness of $\hat A^n$ shows that
$\PP(\Om^n_1)\to0$ as $n\to\iy$. Since $\eps'>0$ is arbitrary,
\eqref{46} follows.

Since rejections occur only when $X^{\circ,n}\ge a^*$, we have
\[
\int 1_{\{X^{\circ,n}(t)<a^*\}}dZ^{\circ,n}(t)=0.
\]
Moreover, we can use \eqref{34} to write
\[
X^{\circ,n}\w a^*=X^{\#,n}(0)+W^{\circ,n}+Y^{\circ,n}-Z^{\circ,n}+E^n,
\qquad E^n=(X^{\circ,n}\w a^*)-X^{\circ,n}.
\]
Combining these relations with \eqref{47} shows that the defining
relations of the Skorohod problem, namely \eqref{44}--\eqref{45}, are valid here, implying
\[
(a^*\w X^{\circ,n}, Y^{\circ,n}, Z^{\circ,n})=\Gam_{[0,a^*]}(X^{\#,n}(0)+W^{\circ,n}+E^n).
\]
By \eqref{46}, $E^n\To0$ uniformly on compacts.
Recall that $W^{\circ,n}$ are $C$-tight. If $\tilde W$ denotes a subsequential limit of it,
using the continuity of $\Gam_{[0,a^*]}$ and using
\eqref{46} once again, shows that along the same subsequence,
$(W^{\circ,n}, X^{\circ,n}, Y^{\circ,n}, Z^{\circ,n})$ converges,
and that its limit satisfies \eqref{43}, as claimed.
The Skorohod map maps continuous paths starting in $[0,a^*]$ to continuous paths.
Hence $(\tilde W,\tilde X,\tilde Y,\tilde Z)$
have continuous paths a.s. This proves the claimed $C$-tightness
of these processes.

{\it Step 2. State space collapse.}
The next major step is to show that the multidimensional process $\hat X^n$
lies close to the minimizing curve. More precisely, we will show that, as $n\to\iy$,
\begin{equation}
  \label{35}
  \Del^n(t):=\hat X^n(t)-\gamma^a(X^{\#,n}(t))\To0,
\end{equation}
uniformly on compacts.

Denote by $\calG=\{x\in\calX:\theta\cdot x\le a^*,x=\gamma^a(\theta\cdot x)\}$
the set of points lying on the minimizing curve, and recall the set
$\pl^+\calX=\{x\in\calX:x_i=b_i \text{ for some } i\}$ corresponding to the
buffer limit boundary. These two compact sets
do not intersect. As a result, there exists $\eps_0>0$
such that for any $0<\eps'<\eps_0$, $\calG^{\eps'}$
and $(\pl^+\calX)^{\eps'}$ do not intersect, where for a set $A\in\R^I$ we denote
\[
A^{\eps'}=\{x:\dist(x,A)\le\eps'\}.
\]
In what follows, it is always assumed that $\eps'<\eps_0$.
Forced rejections occur only at times when $\hat X^n$ lies in $(\pl^+\calX)^{\eps'}$
(for all $n$ large). As a result, as long as the process
$\hat X^n$ lies in $\calG^{\eps'}$, no forced rejections occur.
This observation can be used to deduce that $\sig^n\le\tau^n$, where
\[
\sig^n=\hat\zeta^n\w\zeta^n,
\]
\[
\hat\zeta^n=\inf\{t:X^{\#,n}\ge a^*+\eps'\},\qquad
\zeta^n=\inf\{t:\max_{i\le I}|\Del^n_i(t)|\ge\eps'\}.
\]
Note carefully that $\sig^n$ is not precisely given as
$\inf\{t:\hat X^n(t)\notin\calG^{\eps'}\}$, because $X^{\#,n}$ is defined
using $\theta^n$ while $\gamma^a$ and $\calG$ are defined with $\theta$.
However, since $\theta^n\to\theta$ and $\hat X^n$ remains bounded,
the conclusion that $\sig^n\le\tau^n$,
provided that $n$ is sufficiently large, is valid.

We turn to proving \eqref{35}.
It suffices to show that $\PP(\sig^n<T)\to0$, for any small $\eps'>0$ and any $T$.
Fix $\eps'$ and $T$. Thanks to the fact that $\sig^n\le\tau^n$,
\begin{align}\label{60}
\notag
\PP(\sig^n<T)&=\PP(\sig^n<T,\sig^n\le\tau^n)\\
\notag
&\le\PP(\hat\zeta^n\w\zeta^n\le T\w\tau^n)
\\
&\le\PP(\hat\zeta^n\le T\w\tau^n)+\PP(\zeta^n\le T\w\tau^n).
\end{align}
We have established in Step 1 the convergence \eqref{46}, from which it follows that
$\PP(\hat\zeta^n\le T\w\tau^n)\to0$ as $n\to\iy$.
It therefore suffices to prove the following.

\begin{lemma}\label{lem3}
$\PP(\zeta^n\le T\w\tau^n)\to0$ as $n\to\iy$.
\end{lemma}

\proof
On $\zeta^n\le T\w\tau^n$ let
$x^n:=X^{\#,n}(\zeta^n)=X^{\circ,n}(\zeta^n)$ and let $j=j^n$ and $\xi^n$ be the corresponding
components from the representation $(j,\xi)$ of $x^n$ (with $w=x^n$).

Fix a positive integer $K=K(\eps')=[c_0/\eps']$, where $c_0$ is a constant depending only
on $\theta$, whose value will be specified at a later stage of the proof.
Consider the covering of $[0,\bx]$ by the
$K-1$ intervals $\X_k=\bB(k\eps_1,\eps_1)$, $k=1,2,\ldots,K-1$, where $\bB(x,a)$ denotes
$[x-a,x+a]$ and $\eps_1=\bx/K$. Let also $\tilde\X_k=\bB(k\eps_1,2\eps_1)$.

Recall that $X^{\circ,n}$ are $C$-tight. Invoking the characterization of $C$-tightness
\eqref{105}--\eqref{106}, given $\del>0$ there exists
$\del'=\del'(\del,T,\eps_1)>0$
such that for all sufficiently large $n$,
\begin{equation}
  \label{39}
  |X^{\circ,n}(s)-X^{\circ,n}(t)|\le\eps_1 \text{ for all } s,t\in[0,T], |s-t|\le\del',
  \text{ with probability at least } 1-\del.
\end{equation}
Fix such $\del$ and $\del'$. Then for all large $n$,
\begin{equation}\label{63}
\PP(\zeta^n\le T\w\tau^n)\le\del+\sum_k\PP(\Om^{n,k}),
\end{equation}
where, denoting by $\bT^n$ the interval $[(\zeta^n-\del'\vee 0),\zeta^n]$,
\[
\Om^{n,k}=\{\zeta^n\le T\w\tau^n,x^n\in\X_k,X^{\#,n}(t)\in\tilde\X_k\text{ for all }
t\in\bT^n\}.
\]
Note that we have used the identity $X^{\#,n}=X^{\circ,n}$ on $[0,\tau^n]$.
We fix $k$ and analyze $\Om^{n,k}$, by an argument similar to (but somewhat
more complicated than) that used in Step 1 to treat $\Om^n_1$.

The value assigned by the policy to $B^n$ (see \eqref{32}) remains fixed as
$\hat X^n$ varies within any of the intervals $(\hat a_j,\hat a_{j+1})$.
Aiming at showing that $\PP(\Om^{n,k})\to0$ as $n\to\iy$, for each $k$,
we will first consider the case where $\hat X^n$ remains
in one of these intervals during the time window $\bT^n$;
that is,
\\ (I) $\tilde\X_k\subset(0,a^*)$ and for all $j$,
$\hat a_j\notin\tilde\X_k$. Then we consider the cases
\\ (II) $\tilde\X_k\subset(0,a^*)$
but $\hat a_j\in\tilde\X_k$ for some $j\in\{1,2,\ldots,I-1\}$.
\\ (III) $0\in\tilde\X_k$.
\\ (IV) $a^*\in\tilde\X_k$.

There may be additional intervals $\tilde\X_k$, but they
are all subsets of $(a^*,\iy)$ and therefore not important for our purpose.

(I) $\tilde\X_k\subset(0,a^*)$ and for all $j$, $\hat a_j\notin\tilde\X_k$.
Note that this means that all points $x$ in
$\tilde\X_k$ lead to the same $j$ in the representation $(j,\xi)$ of $x$ given by
\eqref{61}.
Note that $j=j(k)$ depends on $k$ only, and in particular does not vary with $n$.
Also, $j=j^n$ under $\Om^{n,k}$.
In what follows, $j=j(k)$.

Fix $i\in\{j+1,\ldots,I\}$ (unless $i=I$).
We estimate the probability that, on $\Om^{n,k}$, $\zeta^n\le T\w\tau^n$ occurs by
having $\Del^n(\zeta^n)\ge\eps'$.
More precisely, note that $\gamma^a_i(x^n)=a_i$ (because $i>j$).
Then we will show that
\begin{equation}\label{41}
\text{for every } \eps''\in(0,\eps'),\quad
\PP(\Om^{n,k}\cap\{\hat X^n_i(\zeta^n)>a_i+\eps''\})\to0\quad \text{ as } n\to\iy.
\end{equation}
Note that $\gamma^a$ is continuous and that $\Del^n(0)\to0$ as $n\to\iy$,
by \eqref{49}. Using the fact that the jumps
of $\hat X^n$ are of size $n^{-1/2}$, on the event indicated in \eqref{41} there must exist
$\eta^n\in[0,\zeta^n]$ with the properties that
\begin{equation}\label{55}
\hat X^n_i(\eta^n)<a_i+\eps''/2,\qquad X^n_i(t)>a_i \text{ for all } t\in[\eta^n,\zeta^n].
\end{equation}
On this event, during the time interval $[\eta^n,\zeta^n]$,
$i$ is always a member of $\calH(\hat X^n)$, and therefore by
\eqref{32}--\eqref{31}, $B^n_i(t)=\rho'_i(\hat X^n(t))>\rho_i+c$, for some constant $c>0$.
Thus by \eqref{18}, $\frac{d}{dt}\hat Y^n_i\le-\frac{\mu^n_i}{\sqrt n}c$.
Moreover, if we define $\hat\eta^n=\eta^n\vee(\zeta^n-\del')$ then for all
$t\in[\hat\eta^n,\zeta^n]$ one has $\hat X^n(t)\in\tilde\X_k\subset(0,a^*)$
and therefore no rejections occur.
Using these facts in \eqref{17}, we have
\begin{equation}\label{40}
\hat X^n_i[\hat\eta^n,\zeta^n]=\hat W^n_i[\hat\eta^n,\zeta^n]
-c\frac{\mu^n_i}{\sqrt n}(\zeta^n-\hat\eta^n).
\end{equation}
Again, fix a sequence $r_n>0$ with $r^n\to0$ and $r^n\sqrt n\to\iy$.
If $\zeta^n-\eta^n<r_n$ and $n$ is sufficiently large then $\hat\eta^n=\eta^n$,
thus by \eqref{41} and the definition of $\eta^n$, $\hat X^n_i[\hat\eta^n,\zeta^n]\ge\eps''/2$.
As a result,
\[
\bar w_T(\hat W^n_i;r_n)\ge\hat W^n_i[\eta^n,\zeta^n]\ge\eps''/2
\]
must hold.
If, on the other hand, $\zeta^n-\eta^n\ge r_n$ then by \eqref{40},
\[
2\|\hat W^n_i\|_T\ge\hat W^n_i[\hat\eta^n,\zeta^n]\ge c\frac{\mu^n_i}{\sqrt n}r_n
\ge cr_n\sqrt n,
\]
for some constant $c>0$.
Hence the probability in \eqref{41} is bounded by
\begin{equation}
\label{53}
\PP(\bar w_T(\hat W^n_i;r_n)\ge\eps''/2)+\PP(2\|\hat W^n_i\|_T\ge cr_n\sqrt n),
\end{equation}
which converges to zero as $n\to\iy$, by $C$-tightness of $\hat W^n$. This proves \eqref{41}.

Next, if we fix $i<j$ (provided $j\ne1$) then
whenever $\hat X^n_i>0$, $i$ is a member of the high priority set $\calH(\hat X^n)$.
Hence the same argument gives
\begin{equation}\label{50}
\text{for every } \eps''\in(0,\eps'),\quad
\PP(\Om^{n,k}\cap\{\hat X^n_i(\zeta^n)>\eps''\})\to0\quad \text{ as } n\to\iy.
\end{equation}

Consider now $j$ itself. We will show, for the case $j<I$,
\begin{equation}\label{52}
\text{for every } \eps''\in(0,\eps'),\quad
\PP(\Om^{n,k}\cap\{\Del^n_j(\zeta^n)>\eps''\})\to0
\quad \text{ as } n\to\iy.
\end{equation}
Suppose that we show (except in the case $j=I$) for any fixed $\eps''$ and all large $n$ that
on the event indicated in \eqref{52},
\begin{equation}
  \label{51}
  \text{$j\in\calH(\hat X^n(t))$ whenever, prior to $\zeta^n$, one has
  $\Del^n_j(t)\in(\eps''/2,\eps'')$.}
\end{equation}
Then we can argue as in the case of $i>j$, with the following modifications.
Let $C^n(t)=\gamma^a_j(X^{\circ,n}(t))$.
Then $\Del^n_j=\hat X^n_j-C^n$, and
similarly to \eqref{55}, there exists $\eta^n\le\zeta^n$ such that
\[
\Del^n_j(\eta^n)<\eps''/2,\qquad \Del^n_j(t)>0 \text{ for all }
t\in[\eta^n,\zeta^n].
\]
Since by \eqref{51} $j$ is high priority during this interval
we will still have identity \eqref{40} valid.
Arguing separately for the cases $\zeta^n-\eta^n<r_n$ and
$\zeta^n-\eta^n\ge r_n$, leads, in analogy to \eqref{53}, to the conclusion that
the probability in \eqref{52} is bounded by
\begin{equation}
\label{54}
\PP(\bar w_T(\hat W^n_i;r_n)+\bar w_T(C^n;r_n)\ge\eps''/2)
+\PP(2\|\hat W^n_i\|_T+2\|C^n\|_T\ge cr_n\sqrt n).
\end{equation}
In addition to the $C$-tightness of $\hat W$, we now invoke that
of $C^n$, which follows from the continuity of $\gamma^a$ and the $C$-tightness
of $X^{\circ,n}$. This shows \eqref{53}.

Now, since
$\theta\cdot\gamma^a(\theta\cdot x)=\theta\cdot x$ for all $x\in\calX$, $\theta^n\to\theta$,
and $\gamma^a$ uniformly continuous and $\calX$ bounded, we have
\begin{equation}
  \label{56}
  q_n:=\sup_{x\in\calX}|\theta\cdot\gamma^a(\theta^n\cdot x)-\theta\cdot x|\to0,\quad\text{as }
  n\to\iy.
\end{equation}
To show that \eqref{51} holds (except in the case $j=I$),
note by \eqref{56} that
\(
|\theta\cdot\hat X^n(t)-\theta\cdot\gamma^a(X^{\circ,n}(t))|\le q_n\to0.
\)
If $\Del^n_j(t)\ge\eps''/2$ then
\[
-\theta_j\eps''/2\ge\sum_{i\ne j}\theta_i(\hat X^n_i-\gamma^a_i)-\|\theta\|q_n
\ge\sum_{i>j}\theta_i(\hat X^n_i-a_i)-\|\theta\|q_n,
\]
where we used $\gamma^a_i=\gamma^a_i(X^{\circ,n})=0$ for $i<j$ and $\gamma^a_i=a_i$ for $i>j$.
For all large $n$, this implies $\hat X^n_i<a_i$ for at least one $i>j$, by which
$j\in\calH(\hat X^n)$.

We can now show that $\PP(\Om^{n,k})\to0$ as $n\to\iy$. Indeed, in the case $j=I$,
we have by \eqref{50}, using $\gamma^a_i=0$,
$\PP(\Om^{n,k}\cap\{\max_{i<I}|\Del^n_i(\zeta^n)|>\eps''\})\to0$.
By \eqref{56}, $|\theta\cdot\Del^n(\zeta^n)|\le q_n$.
Since $\theta\in(0,\iy)^I$ and $q_n\to0$, this shows that
\begin{equation}
\label{57}
\PP(\Om^{n,k}\cap\{\max_{i\le I}|\Del^n_i(\zeta^n)|>\eps''\})\to0.
\end{equation}
In the case $j<I$, combining \eqref{41}, \eqref{50}, \eqref{52},
we have $\PP(\Om^{n,k}\cap\{\max_{i\le I}\Del^n_i(\zeta^n)>\eps''\})\to0$.
Using again the fact $|\theta\cdot\Del^n(\zeta^n)|\le q_n\to0$ gives
that \eqref{57} is valid in this case as well.

Since $\eps''$ is arbitrarily small, it follows from the definition of $\zeta^n$
that $\PP(\Om^{n,k})\to0$ as $n\to\iy$.

(II) $\tilde\X_k\subset(0,a^*)$ but $\hat a_j\in\tilde\X_k$ for some $j\in\{1,2,\ldots,I-1\}$.
Let $(j^n(t),\xi^n(t))$ denote the representation \eqref{61} for $X^{\#,n}(t)$.
The difficulty here is that in the
time window $\bT^n$, $j^n$ varies between two values, namely $j$ and $j+1$,
and it is no longer true that
$\gamma^a_{j+1}(X^{\#,n})=a_{j+1}$ on that time interval.
The way we treat this is by bounding $\Del^n$ from above by a quantity that depends on
$\eps_1$, rather than by an arbitrarily small $\eps''$.
To this end, let us show that on $\Om^{n,k}$,
\begin{equation}\label{64}
\gamma^a_{j+1}(X^{\#,n}(t))\ge
a_{j+1}-c_1\eps_1,\qquad
t\in\bT^n,
\end{equation}
where $c_1=4/\theta_{\rm min}$ and $\theta_{\rm min}=\min_i\theta_i$.
Indeed, we have for any $w\in\tilde\X_k$, $|w-\hat a_j|\le 4\eps_1$, since
$\hat a_j$ is also in $\tilde\X_k$. Now, if $w\ge\hat a_j$
then $\gamma^a_{j+1}(w)=a_{j+1}$. Otherwise,
\[
w=\hat a_{j+1}+\theta_{j+1}\xi=\hat a_j-\theta_{j+1}a_{j+1}+\theta_{j+1}\xi,
\]
thus $|a_{j+1}-\xi|\le4\theta_{j+1}^{-1}\eps_1$, whence follows \eqref{64}.

Now, \eqref{41} is valid for all $i>j+1$, by the proof given in case (I).
For $i=j+1$ it is also valid, even though $\gamma^a_{j+1}(X^{\#,n}(t))$ is not
necessarily equal to $a_{j+1}$.
For $i<j$, \eqref{50} is valid with the same proof. As for $i=j$, \eqref{52}
is valid with same proof
(the fact that $\gamma^a_j$ may assume the value zero does not affect this proof).

Combining all the estimates except for $i=j+1$ gives, for all small $\eps''$,
\[
\PP(\Om^{n,k}\cap\{\max_{i\ne j+1}\Del^n_i(\zeta^n)>\eps''\})\to0.
\]
For $i=j+1$, the estimate \eqref{41} and the bound \eqref{64} give
$\PP(\Om^{n,k}\cap\{\Del^n_{j+1}(\zeta^n)>2c_1\eps_1\})\to0$ as $n\to\iy$.
Along with $|\theta\cdot\Del^n(\zeta^n)|\le q_n$, this gives
\[
\PP(\Om^{n,k}\cap\{\max_{i\le I}|\Del^n_i(\zeta^n)|>3c_1\eps_1\})\to0,
\]
as $n\to\iy$. We now determine the constant $c_0$ used to define $K$.
We do so in such a way that $3c_1\eps_1<\eps'/2$. In particular, any constant $c_0>6c_1\bx
=24\bx/\theta_{\rm min}$ will do.
This way we obtain $\PP(\Om^{n,k})\to0$ as $n\to\iy$.

(III) $0\in\tilde\X_k$.
The only difference of this case from case (I) is that during $\bT^n$, $\hat X^n$ may hit
zero, and so by \eqref{58} and \eqref{32}, $B^n$ will be zero. However, the analysis in
case (I) is performed only on intervals where $\hat X^n\ne0$, and as a result gives
rise to the same conclusion, namely
$\PP(\Om^{n,k})\to0$ as $n\to\iy$.

(IV) $a^*\in\tilde\X_k$. In this case, during $\bT^n$, $\theta\cdot\hat X^n$ may exceed
$a^*$, and so rejections of class $i^*$ customers may occur.
The only way it affects the proof of case (I)
is by adding a negative term to the r.h.s.\ of \eqref{40}. However, the consequences
of \eqref{40} remain valid with this addition.
(Note that
for all sufficiently small $\eps$ one has $\hat a_i\ne a^*$ for all $i$, hence assuming
$\eps$ is sufficiently small, we do not need to check case (II) here.)

Having shown that $\PP(\Om^{n,k})\to0$ in all cases, using \eqref{63} and the fact that
$\del>0$ is arbitrary completes the proof of the lemma.
\qed

As a consequence of the lemma and \eqref{60}, we have $\PP(\sig^n<T)\to0$ as
$n\to\iy$. Since $\eps'$ is arbitrary, \eqref{35} is established.

{\it Step 3. Weak convergence.}
Having shown that $\PP(\sig^n<T)\to0$, we have, using $\sig^n\le \tau^n$,
that $\PP(\tau^n<T)\to0$.
As a result, the conclusion of Step 1 regarding the stopped processes
holds also for the unstopped ones. That is,
$(W^{\#,n},X^{\#,n},Y^{\#,n},Z^{\#,n})$
are $C$-tight, and any subsequential limit satisfies \eqref{43} a.s.

Let $W$ be an $(m,\sig)$-BM (of dimension $I$) and set $\bar W=\theta\cdot W$.
Denote by $(\bar X,\bar Y,\bar Z)$ the triple from Proposition \ref{prop2},
i.e., $(\bar X,\bar Y,\bar Z)=\Gam_{[0,\bx^*]}(\bar x_0+\bar W)$ (note that $\bar W$
is a $(\bar m,\bar\sigma)$-BM).
Also let $(\bar X^a,\bar Y^a,\bar Z^a)=\Gam_{[0,a^*]}(\bar x_0+\bar W)$.

For any finite $T$, the sequence $Z^{\#,n}(T)$ is tight. On the event $\tau^n>T$,
which has overwhelming probability,
\begin{equation}\label{65}
\hat Z^n(T)=\hat Z^n_{i^*}(T)e^{(i^*)},
\end{equation}
hence $\|\hat Z^n(T)\|$ is a tight sequence. The bound \eqref{30} thus gives the
tightness of $\|\hat Y^n(T)\|$.
The argument from the lower bound in the paragraph following \eqref{30}
shows that $\hat W^n\To W$ as $n\to\iy$.
Thus \eqref{43} determines the limit of the one-dimensional processes,
namely $(W^{\#,n},X^{\#,n},Y^{\#,n},Z^{\#,n})\To(\bar W,\bar X^a,\bar Y^a,\bar Z^a)$.
Moreover,
$(\hat W^n,\hat X^n,\hat Y^n,\hat Z^n)\To(W,X,Y,Z)$ where
$\theta\cdot X=\bar X^a$ and $\gamma^a(\bar X^a)$ (by \eqref{35})
$Z=\zeta^*\bar Z^a$ (by \eqref{65}) and $Y=X-x_0-W+Z$, by \eqref{17}.
We obtain precisely
the relations from Proposition \ref{prop1}, except that
the reflection interval is $[0,a^*]$ rather than $[0,\bx^*]$.

We have shown that, as $n\to\iy$,
\[
\int_0^\iy e^{-\al t}[h\cdot\hat X^n(t)+\al r\cdot \hat Z^n(t)]dt]\To
\int_0^\iy e^{-\al t}[h\cdot\gamma^a(\bar X^a(t))+\al \bar r \bar Z^a(t)]dt.
\]

{\it Step 4. Convergence of costs.}
Since $\hat X^n$ are uniformly bounded, we immediately obtain
$\EE\int_0^\iy e^{-\al t}h\cdot\hat X^n(t)dt\to
\EE\int_0^\iy e^{-\al t}h\cdot\gamma^a(\bar X^a(t))dt$.
As for the second term, we borrow an argument from \cite{bell-will-1}.
Consider the probability space
$(\R_+\times\Om,\calB(\R_+)\times\calF,m\times\PP)$,
where $dm=\al e^{-\al t}dt$.
Then the result of the previous
step can be expressed as the convergence in law,
\[
r\cdot\hat Z^n\to \bar r \bar Z^a,
\]
w.r.t.\ the probability measure $m\times\PP$.
Thus to obtain $\EE\int_0^\iy e^{-\al t}r\cdot\hat Z^n(t)dt
\to \EE\int_0^\iy e^{-\al t}\bar r \bar Z^a(t)dt$,
it suffices to show the $m\times\PP$-uniform integrability (UI) of $r\cdot\hat Z^n$.
For this, it suffices that
\begin{equation}\label{66}
\limsup_n\EE\int_0^\iy e^{-\al t}\|\hat Z^n(t)\|^2dt<\iy.
\end{equation}
It is established in equation (172) of \cite{bell-will-1} that
\begin{equation}\label{67}
\EE[(\|\hat S^n\|_t)^2]\le c(1+t)
\end{equation}
for a constant $c$ independent of $n$ and $t$,
with the same estimate holding for $\hat A^n$. In what follows, we show that we can deduce
\eqref{66} from \eqref{67}.

To this end, recall that rejections
occur only when either $\theta\cdot\hat X^n\ge a^*$ or, for some $i$,
$\hat X^n_i\ge a_i-n^{-1/2}$.
In particular, if we let $\bar a= a^*\w\min_i(a_i\theta_i/2)$,
then using the convergence $\theta^n\to\theta$, we have, for all large $n$,
that no rejections take place when $X^{\#,n}=\theta^n\cdot\hat X^n<\bar a$.
Consider the truncated version $X^{1,n}:=\bar a\w X^{\#,n}$ of $X^{\#,n}$. Then by \eqref{34},
\begin{equation}\label{59}
X^{1,n}(t)=W^{1,n}+Y^{\#,n}-Z^{\#,n},
\end{equation}
where we denote
\[
W^{1,n}=W^{\#,n}+E^n,\qquad E^n=X^{\#,n}(0)+X^{1,n}-X^{\#,n}.
\]
By the above discussion,
\[
\int_0^\iy 1_{\{X^{1,n}<\bar a\}}dZ^{\#,n}=0,
\]
and by the work conservation property,
\[
\int_0^\iy 1_{\{X^{1,n}>0\}}dY^{\#,n}=0.
\]
Moreover, the initial value $W^{1,n}(0)$ lies in $[0,\bar a]$.
These facts dictate that $(X^{1,n},Y^{\#,n},Z^{\#,n})$ solves
the SP on $[0,\bar a]$ for $W^{1,n}$. That is,
\(
(X^{1,n},Y^{\#,n},Z^{\#,n})=\Gam_{[0,\bar a]}(W^{1,n}).
\)
It is well-known (see e.g., \cite{KLRS}) that $\Gam_{[0,\bar a]}$ is uniformly Lipschitz
in the following strong sense: There exists a constant $L$ depending only on $\bar a$,
such that for every
$w_1,w_2\in D([0,t],\R)$, one has
$\|x_1-x_2\|_t+\|y_1-y_2\|_t+\|z_1-z_2\|_t\le L\|w_1-w_2\|_t$,
where $(x_i,y_i,z_i)=\Gam_{[0,\bar a]}(w_i)$, $i=1,2$. Since the response to $0$ is
$(0,0,0)$, it follows that
\[
Z^{\#,n}(t)\le L\|W^{1,n}\|_t\le c(\|W^{\#,n}\|_t+2\bx),
\qquad t\ge0,
\]
where we used the bound $|E^n|\le 2\bx$.
Since $\theta^n$ converge to $\theta\in(0,\iy)^I$, this and the definitions
\eqref{33} imply
\[
\|\hat Z^n(t)\|\le c(1+\|\hat W^n\|_t).
\]
Going back to \eqref{26} and using the fact that $T^n_i(t)\le t$ for each $i$,
and the convergence of $m^n_i$ (see \eqref{19}),
\[
\|\hat Z^n(t)\|\le c(1+t+\|\hat A^n\|_t+\|\hat S^n\|_t),
\qquad t\ge0,
\]
where $c$ is independent of $n$ and $t$. Combining this with \eqref{67} gives \eqref{66}.
Hence follows the required UI.

We have thus proved that, with $U^n=U^n(\eps)$,
\[
\lim_nJ^n(U^n)=\lim_n\E\Big[\int_0^\iy e^{-\al t}[h\cdot\hat X^n(t)+\al r\cdot\hat Z^n(t)]dt\Big]
=\E\Big[\int_0^\iy e^{-\al t}[h\cdot\gamma^a(\bar X^a(t))+\bar r\bar Z^a(t)]dt\Big].
\]
Denoting the right member above by $V(x_0;\eps)$, the result will follow once we prove
that $V(x_0,\eps)\to V(x_0)$ as $\eps\to0$.
Now, as $\eps\to0$, one has $a\to b$, $a^*\to\bx^*$ and $\gamma^a\to\gamma$ uniformly.
Moreover, the process $(\bar X^a,\bar Y^a,\bar Z^a)$
converges to $(\bar X,\bar Y,\bar Z)$ in law, as can be deduced, for example, from
the explicit representation of $\Gamma|_{[0,a]}$ provided in \cite{KLRS}.
Thus
\[
h\cdot\gamma^a(\bar X^a)+\al\bar r\bar Z^a\to h\cdot\gamma(\bar X)+\al\bar r\bar Z
=\bar h(\bar X)+\al\bar r\bar Z,
\]
in law w.r.t.\ $m\times\PP$.
Now, $\gamma^a$ is bounded; hence to prove the convergence
\[
V(x_0,\eps)\to\check V(x_0):=
\E\int_0^\iy e^{-\al t}(\bar h(\bar X)+\al\bar r\bar Z)dt,
\]
it suffices to show the corresponding UI, and in particular, that
\begin{equation}\label{68}
\limsup_{a^*\to\bx^*}\EE\int_0^\iy e^{-\al t}\|\bar Z^a(t)\|^2dt<\iy.
\end{equation}
To see that \eqref{68} holds,
apply Ito's formula to $(\bar X^a(t))^2$, use the facts
$\int_0^t\bar X^a(s)d\bar Y^a(s)=0$
and $\int_0^t\bar X^a(s)d\bar Z^a(s)=aZ^a(t)$, to get
\[
Z^a(t)=\frac{1}{2a}\Big\{(\bar X^a(0))^2-(\bar X^a(t))^2+2\int_0^t\bar X^a(s)d\bar W(s)
+\bar\sigma^2t\Big\}.
\]
Since $\bar X^a$ is bounded by $a$, \eqref{68} follows easily.
As a result we have $V(x_0,\eps)\to\check V(x_0)$ as $\eps\to0$.
Using integration by parts, $\check V(x_0)=\E\int_0^\iy e^{-\al t}(\bar h(\bar X)dt
+\bar r d\bar Z(t))$. According to Proposition \ref{prop2}, this is precisely $\bar V(\bar x_0)$,
because $(\bar Y,\bar Z)$ is optimal for $\bar V(\bar x_0)$.
By Proposition \ref{prop1}, $\bar V(\bar x_0)=V(x_0)$.
This proves the statement of the theorem.

{\it Step 5. General initial condition.}
Finally, we relax the assumption \eqref{49} on the initial condition.
Here we do not give the proof in full detail, but only a brief sketch.
Let $\eps$ be given and let $a$, $a^*$ be as before.
Let $\tau^n_0$ denote the first time when a condition analogous to \eqref{49} holds; more precisely,
let $\alpha_n>0$, $\alpha_n\to0$ and
\[
\tau^n_0=\inf\{t:\|\hat X^n(t)-\gamma^a(\theta\cdot\hat X^n(t))\|\le\alpha_n
\text{ and } \theta^n\cdot\hat X^n(t)\in[0,a^*]\}.
\]
The idea is to show that
(i) with a suitable choice of $\alpha_n$, one has $\tau^n_0\to0$ in probability, and
(ii) starting from $(\tau^n_0,\hat X^n(\tau^n))$ in place of $(0,\hat X^n(0))$,
the arguments in all the proof can be repeated without additional effort.
While (ii) follows in a straightforward manner, but notationally heavy, (i) is
a consequence of similar to the proof of \eqref{35}. We omit the details.
\qed

\section{On pathwise Little's law and throughput time constraints}\label{sec5}
\beginsec

This section is motivated by the the work of
Plambeck et al.\ \cite{PKH}, where cumulative rejection costs are minimized
subject to throughput time constraints in heavy traffic.
The pathwise Little's law (see \cite{pw-little} for its original version)
and the related Reiman's snapshot principle \cite{rei-snap}
imply that throughput time, regarded as a process,
multiplied by the arrival rate is asymptotically equal to the
queuelength process in the heavy traffic limit
(Proposition \ref{prop3} below provides a precise statement for the present model).
It is thus natural to expect that the solution to the problem with finite
buffers that we have given
can be transformed into one where throughput times are constrained, and
rejection and holding costs are incurred.
Our results in this direction are only partial. The main purpose of this section
is to propose this problem setting and comment on relations to the main body
of the paper, leaving the main question open.

Some additional notation is necessary in order to formulate the
throughput time constraint problem.
Recall that $A^n_i$ and $Z^n_i$ denote the arrival and rejection counting processes,
and let $AD^n_i=A^n_i-Z^n_i$ denote the admission counting processes. Given $t\ge0$ and
$i\in\calI$ we denote by $AT^n_i(t)$ (where $AT$ is mnemonic for arrival time)
the first time after $t$ when a customer of class $i$
arrives and is admitted into the system, namely
\[
AT^n_i(t)=\inf\{s>t: AD^n_i(s)>AD^n_i(t)\},\qquad t\ge0.
\]
Let $DT^n_i(t)$ denote the departure time of that customer.
This process can be recovered from the other processes we have defined, as follows
\begin{equation}\label{500}
DT^n_i(t)=\inf\{s:D^n_i(s)\ge D^n_i(AT^n_i(t))+X^n_i(AT^n_i(t))\}.
\end{equation}
To see this, note that the time of departure of a given customer
equals the time when all customers of its class present in the system at the time of its arrival
(including the given customer) have departed the system.
This gives the identity
\begin{equation}\label{501}
D^n_i(DT^n_i(t))-D^n_i(AT^n_i(t))=D^n_i[AT^n_i(t),DT^n_i(t)]=X^n_i(AT^n_i(t)),
\end{equation}
by which \eqref{500} follows. (Note, by right continuity, that indeed
$X^n_i(AT^n_i(t))$ equals the number of class-$i$ customers in the system at the time
of arrival $AT^n_i(t)$, which includes that customer).

The time the customer spends in the system, that we will call
the {\it throughput time}, is given by
\begin{equation}\label{502}
\Th^n_i(t)=DT^n_i(t)-AT^n_i(t).
\end{equation}
Denote a diffusion scale version of the throughput time by
\[
\hat\Th^n_i(t)=\sqrt n\Th^n_i(t).
\]
Recall the distinction between an admissible control and an admissible
control satisfying the buffer constraints (Definition \ref{def1}) and the corresponding
classes $\tilde\calU^n$ and $\calU^n$. In this section we replace the buffer constraint
by an asymptotic requirement on the throughput times.
Following \cite{PKH}, we fix constants $d_i>0$, $i\in\calI$, and introduce
\begin{definition}\label{def4} {\bf (Asymptotic compliance)}
A sequence $\{U^n\}_{n\in\N}$, $U^n\in\tilde\calU^n$, of admissible controls, satisfying
\begin{equation}
  \label{70}
  \text{for every $T$, } \max_i\sup_{t\in[0,T]}(\hat{\Th}^n_i(t)-d_i)^+\To0,
\end{equation}
is said to be {\bf asymptotically compliant}. Denote by $\cal AC$ the set of all
asymptotically compliant sequences of admissible controls.
\end{definition}
Plambeck et al. \cite{PKH} study AO
in presence of rejection costs. In what seems to be a natural extension
of their problem to include holding costs, we consider
minimizing $J^n$ among all policies satisfying \eqref{70} instead of the buffer
constraints. Thus we set
\[
V^{\cal AC}=\inf_{\{U^n\}\in\cal AC}\liminf_{n\to\iy}J^n(U^n).
\]
We address this problem by comparing it to the problem with buffer constraints
via a conditional pathwise Little's law, presented next.
Its proof is deferred to the end of the section.
\begin{proposition}\label{prop3}
Fix $T>0$.
Given any sequence of admissible controls $U^n\in\tilde\calU^n$, if
$\hat Y^n$ are $C$-tight
and $\|\hat X^n\|_T$ are tight then $E^n\To0$ uniformly over $[0,T]$, where
\[
E^n_i(t)=\hat X^n_i(AT^n_i(t))-\la_i \hat\Th^n_i(t).
\]
If, in addition, $\hat X^n$ are $C$-tight then also $\tilde E^{n}\To0$, where
\begin{equation}\label{92}
\tilde E^n_i(t)=\hat X^n_i(t)-\la_i \hat\Th^n_i(t).
\end{equation}
\end{proposition}
This suggests that a constraint $d_i$ on $\hat\theta^n_i$ should be similar
a constraint $b_i=\la_id_i$ on $\hat X^n$, as in the formulation with finite buffers.
Recall that $V$ denotes the value function for the RBCP.
This definition depends upon the choice
of the set $\calX$, which we take to be rectangular as in \eqref{10}, with $b_i=\la_id_i$.
Recall the policies defined before, and in particular, Remark \ref{rem3} by which
the policy with $\eps=0$ makes perfect sense. Here we do not have
strict buffer constraints, only the requirement to meet the asymptotic
constraint \eqref{70}, therefore we can work with simply $\eps=0$.
We denote the resulting policy by $U^n_*$.

\begin{theorem}
  \label{th3}
The sequence $\{U^n_*\}$ is asymptotically compliant.
  Moreover, $\limsup_{n\to\iy}J^n(U^n_*)\le V(x_0)$.
\end{theorem}

\proof
Step 4 in the proof of Theorem \ref{th2} gives the upper bound
on the cost. Thus it suffices to show asymptotic compliance.
Now, it follows from the proof of Theorem \ref{th2} that the controls under
consideration satisfy the assumptions of Proposition \ref{prop3}.
Using this proposition along with the fact that $\PP(\sigma^n<T)\to0$
(see the proof of Theorem \ref{th2}), gives the result.
\qed

\begin{conjecture}
  \label{conj1}
  One has $V^{\cal AC}\ge V(x_0)$.
\end{conjecture}

If the above is true then, by Theorem \ref{th3},
$\{U^n_*\}$ are AO for the problem under consideration.
One might approach the conjecture by using Proposition \ref{prop3}
to connect to the lower bound of Theorem \ref{th1}.
The difficulty here is that one must consider an arbitrary sequence
of controls, and there is no guarantee that the assumptions of Proposition \ref{prop3},
particularly the $C$-tightness of $\hat Y^n$, hold in such generality.
We are able to show a partial result.

We address only policies that give rise to state space collapse.
More precisely, consider a sequence $\{U^n\}_{n\in\N}\in{\cal AC}$, and write $\{U^n\}\in\widehat{\cal AC}$
if it satisfies the following. (i) Each $U^n$ is
work conserving;
(ii) for some $\hat\bx\in[0,\hat\bx)$, rejections occur only when the scaled workload exceeds $\hat\bx$,
and only from one particular class (save forced rejections); and
(iii) for some continuous $\hat\gamma:[0,\bx]\to\calX$ satisfying
\begin{equation}\label{91}
\{x\in\calX:\theta\cdot x\le\hat\bx,x=\hat\gam(\theta\cdot x)\}\cap
\pl^+\calX=\emptyset,
\end{equation}
one has $\hat X^n-\hat\gam(\theta\cdot\hat X^n)\To0$ as $n\to\iy$.
Set
\[
V^{\widehat{\cal AC}}=\inf_{\{U^n\}\in\widehat{\cal AC}}\liminf_{n\to\iy} J^n(U^n).
\]
\begin{proposition}
  \label{prop4}
  One has $V^{\widehat{\cal AC}}\ge V(x_0)$.
\end{proposition}
This result is far from being satisfactory. However, it shows that the two
formulations are equivalent at least for this restricted class of policies.
The proof is based on various elements of the proofs of the finite buffer results.

\proof
Assume, without loss of generality, that $V^{\widehat{\cal AC}}<\iy$, and consider
$\{U^n\}\in\widehat{\cal AC}$ with $\liminf J^n(U^n)<\iy$. Finiteness of this quantity gives,
along the lines of the proof of Theorem \ref{th1}, that $\hat W^n$ are $C$-tight.
The assumptions on $\{U^n\}$ as a sequence in $\widehat{\cal AC}$ imply, by arguments as in the proof of
Theorem \ref{th2}, that the one-dimensional processes
$(W^{\circ,n},X^{\circ,n},Y^{\circ,n},Z^{\circ,n})$ are
$C$-tight, and any subsequential limit $(\tilde W,\tilde X,\tilde Y,\tilde Z)$
satisfies a.s., \eqref{43}. The state space collapse assumption, along with \eqref{91}
imply that, for any $T$, $\PP(\tau^n<T)\to0$, and that the un-stopped processes
$(X^{\#,n},Y^{\#,n},Z^{\#,n})$ as well as $\hat X^n$ are $C$-tight. Finally,
$C$-tightness of $\hat Z^n$ follows from that of $Z^{\#,n}$ by arguments as in the same proof,
and that of the processes $\hat Y^n$ follows equation \eqref{17} now that we have $C$-tightness
of all the other processes involved. This verifies the assumptions of Proposition \ref{prop3}.

As a result $\tilde E^n\To0$ (where $\tilde E^n$ are defined in \eqref{92}).
The asymptotic compliance of the sequence of controls along
with the convergence $\tilde E^n\To0$ imply the validity of the
relaxed assumption \eqref{110} under which the lower bound, Theorem \ref{th1}, is proved. Thus we conclude
from Theorem \ref{th1} that $\liminf J^n(U^n)\ge V(x_0)$ for any $\{U^n\}\in\widehat{\cal AC}$.
\qed

\noi{\bf Proof of Proposition \ref{prop3}.}
This proof is close to that of Lemma A.4 in the e-companion \cite{ata-nds-app}
of \cite{ata-nds}.
By \eqref{501} and \eqref{502},
\(
\hat X^n_i(AT^n_i(t))=n^{-1/2}D^n_i[AT^n_i(t),AT^n_i(t)+\Th^n_i(t)].
\)
Now, by \eqref{1} and \eqref{16-}
\[
D^n_i(t)=S^n_i(T^n_i(t)),\qquad
\tilde V^n_i(t):=\hat S^n_i(T^n_i(t))=\frac{S^n_i(T^n_i(t))-\mu^n_iT^n_i(t)}{\sqrt n},
\]
and
\(
n^{-1/2}{D^n_i}=\sqrt n \bar\mu^n_i T^n_i+\tilde V^n_i.
\)
By \eqref{18}, recalling that $\bar\mu^n_i=\mu^n_i/n$, we have
\[
\sqrt n \bar\mu^n_i T^n_i(t)=\sqrt n\bar\mu^n_i\rho_it-\hat Y^n_i(t).
\]
Hence
\[
\frac{D^n_i(t)}{\sqrt n}= \sqrt n\bar\mu^n_i\rho_it-\hat Y^n_i(t) +\tilde V^n_i(t),
\]
and therefore
\begin{equation}\label{300}
\hat X^n_i(AT^n_i(t))=\bar\mu^n_i\rho_i\hat\Th^n_i(t)
+e^{1,n}_i(t)-e^{2,n}_i(t),
\end{equation}
where
\[
e_i^{1,n}(t)=\tilde V^n_i[AT^n_i(t),AT^n_i(t)+\Th^n_i(t)],\qquad
e_i^{2,n}(t)=\hat Y^n_i[AT^n_i(t),AT^n_i(t)+\Th^n_i(t)].
\]
Now, $\tilde V^n$ is $C$-tight by $C$-tightness of $\hat S^n_i$ and the uniform Lipschitz property
of $T^n_i$. The processes $\hat Y^n$ are assumed to be $C$-tight. Thus by the assumption on
$\|\hat X^n\|_T$, it follows that $\|\hat\Th^n\|_T=\sqrt n\|\Th^n\|_T$ are tight r.v.s,
and thus that $e^{1,n}_i$ and $e^{2,n}_i$ all converge to zero uniformly.
Since $\bar\mu^n_i\rho_i\to\la_i$, this shows the first statement of the result.
The second statement now follows by the uniform convergence of $AT^n_i(t)\to t$.
\qed

\noi{\bf Acknowledgement.} We would like to thank Haya Kaspi for referring us
to the argument in \cite{DM-english} and the referees for helpful comments.

\footnotesize

\bibliographystyle{is-abbrv}


\begin{thebibliography}{10}

\bibitem{ata06}
B.~Ata.
\newblock Dynamic control of a multiclass queue with thin arrival streams.
\newblock {\em Operations Research}, 54\penalty0 (5):\penalty0 876--892, 2006.

\bibitem{AK1}
B.~Ata and S.~Kumar.
\newblock Heavy traffic analysis of open processing networks with complete
  resource pooling: asymptotic optimality of discrete review policies.
\newblock {\em Ann. Appl. Probab.}, 15\penalty0 (1A):\penalty0 331--391, 2005.

\bibitem{ata-ols}
B.~Ata and T.~L. Olsen.
\newblock Near-optimal dynamic lead-time quotation and scheduling under
  convex-concave customer delay costs.
\newblock {\em Operations Research}, 57\penalty0 (3):\penalty0 753--768, 2009.

\bibitem{ata-nds}
R.~Atar.
\newblock A diffusion regime with nondegenerate slowdown.
\newblock {\em Oper. Res.}, 60\penalty0 (2):\penalty0 490--500, 2012.

\bibitem{ata-nds-app}
R.~Atar.
\newblock A diffusion regime with nondegenerate slowdown: Appendix.
\newblock {\em Electronic companion, downloadable at
  http://webee.technion.ac.il/people/atar/online-appendix.pdf}, 2012.

\bibitem{AB2006}
R.~Atar and A.~Budhiraja.
\newblock Singular control with state constraints on unbounded domain.
\newblock {\em Ann. Probab.}, 34\penalty0 (5):\penalty0 1864--1909, 2006.

\bibitem{ABW}
R.~Atar, A.~Budhiraja, and R.~J. Williams.
\newblock H{JB} equations for certain singularly controlled diffusions.
\newblock {\em Ann. Appl. Probab.}, 17\penalty0 (5-6):\penalty0 1745--1776,
  2007.

\bibitem{atagur}
R.~Atar and I.~Gurvich.
\newblock Scheduling parallel servers in the non-degenerate slowdown diffusion
  regime: Asymptotic optimality results.
\newblock {\em Ann. Appl. Probab., to appear}, 2013.

\bibitem{AMR}
R.~Atar, A.~Mandelbaum, and M.~I. Reiman.
\newblock Scheduling a multi class queue with many exponential servers:
  asymptotic optimality in heavy traffic.
\newblock {\em Ann. Appl. Probab.}, 14\penalty0 (3):\penalty0 1084--1134, 2004.

\bibitem{atasol}
R.~Atar and N.~Solomon.
\newblock Asymptotically optimal interruptible service policies for scheduling
  jobs in a diffusion regime with nondegenerate slowdown.
\newblock {\em Queueing Systems Theory Appl.}, 69\penalty0 (217--235), 2011.

\bibitem{bell-will-1}
S.~L. Bell and R.~J. Williams.
\newblock Dynamic scheduling of a system with two parallel servers in heavy
  traffic with resource pooling: asymptotic optimality of a threshold policy.
\newblock {\em Ann. Appl. Probab.}, 11\penalty0 (3):\penalty0 608--649, 2001.

\bibitem{ben-men}
M.~N. Bennani and D.~A. Menasce.
\newblock Resource allocation for autonomic data centers using analytic
  performance models.
\newblock In {\em Autonomic Computing, 2005. ICAC 2005. Proceedings. Second
  International Conference on}, pages 229--240. IEEE, 2005.

\bibitem{Bill}
P.~Billingsley.
\newblock {\em Convergence of probability measures}.
\newblock Wiley Series in Probability and Statistics: Probability and
  Statistics. John Wiley \& Sons Inc., New York, second edition, 1999.
\newblock ISBN 0-471-19745-9.
\newblock x+277 pp.
\newblock A Wiley-Interscience Publication.

\bibitem{bramson-ssc}
M.~Bramson.
\newblock State space collapse with application to heavy traffic limits for
  multiclass queueing networks.
\newblock {\em Queueing Systems Theory Appl.}, 30\penalty0 (1-2):\penalty0
  89--148, 1998.

\bibitem{BG2006}
A.~Budhiraja and A.~P. Ghosh.
\newblock Diffusion approximations for controlled stochastic networks: an
  asymptotic bound for the value function.
\newblock {\em Ann. Appl. Probab.}, 16\penalty0 (4):\penalty0 1962--2006, 2006.

\bibitem{BG2012}
A.~Budhiraja and A.~P. Ghosh.
\newblock Controlled stochastic networks in heavy traffic: convergence of value
  functions.
\newblock {\em Ann. Appl. Probab.}, 22\penalty0 (2):\penalty0 734--791, 2012.

\bibitem{BVW}
C.~Buyukkoc, P.~Varaiya, and J.~Walrand.
\newblock The {$c\mu$} rule revisited.
\newblock {\em Adv. in Appl. Probab.}, 17\penalty0 (1):\penalty0 237--238,
  1985.

\bibitem{cox-smi}
D.~R. Cox and W.~L. Smith.
\newblock {\em Queues}.
\newblock Methuen's Monographs on Statistical Subjects. Methuen \& Co. Ltd.,
  London; John Wiley \& Sons Inc., New York, 1961.
\newblock xii+180 pp.

\bibitem{daidai}
J.~G. Dai and W.~Dai.
\newblock A heavy traffic limit theorem for a class of open queueing networks
  with finite buffers.
\newblock {\em Queueing Systems Theory Appl.}, 32\penalty0 (1-3):\penalty0
  5--40, 1999.

\bibitem{DM-english}
C.~Dellacherie and P.-A. Meyer.
\newblock {\em Probabilities and potential}, volume~29 of {\em North-Holland
  Mathematics Studies}.
\newblock North-Holland Publishing Co., Amsterdam-New York; North-Holland
  Publishing Co., Amsterdam-New York, 1978.
\newblock ISBN 0-7204-0701-X.
\newblock viii+189 pp.

\bibitem{durrett-book}
R.~Durrett.
\newblock {\em Probability: theory and examples}.
\newblock Duxbury Press, Belmont, CA, second edition, 1996.
\newblock ISBN 0-534-24318-5.
\newblock xiii+503 pp.

\bibitem{ethkur}
S.~N. Ethier and T.~G. Kurtz.
\newblock {\em Markov processes}.
\newblock Wiley Series in Probability and Mathematical Statistics: Probability
  and Mathematical Statistics. John Wiley \& Sons Inc., New York, 1986.
\newblock ISBN 0-471-08186-8.
\newblock x+534 pp.
\newblock Characterization and convergence.

\bibitem{gha-war}
S.~Ghamami and A.~R. Ward.
\newblock Dynamic scheduling of a two-server parallel server system with
  complete resource pooling and reneging in heavy traffic: Asymptotic
  optimality of a two-threshold policy.
\newblock {\em Mathematics of Operations Research}, 38\penalty0 (4):\penalty0
  761--824, 2013.

\bibitem{GW1}
A.~P. Ghosh and A.~P. Weerasinghe.
\newblock Optimal buffer size for a stochastic processing network in heavy
  traffic.
\newblock {\em Queueing Syst.}, 55\penalty0 (3):\penalty0 147--159, 2007.

\bibitem{har1988}
J.~M. Harrison.
\newblock Brownian models of queueing networks with heterogeneous customer
  populations.
\newblock In {\em Stochastic differential systems, stochastic control theory
  and applications ({M}inneapolis, {M}inn., 1986)}, volume~10 of {\em IMA Vol.
  Math. Appl.}, pages 147--186. Springer, New York, 1988.

\bibitem{har2000}
J.~M. Harrison.
\newblock Brownian models of open processing networks: canonical representation
  of workload.
\newblock {\em Ann. Appl. Probab.}, 10\penalty0 (1):\penalty0 75--103, 2000.

\bibitem{har2003}
J.~M. Harrison.
\newblock A broader view of {B}rownian networks.
\newblock {\em Ann. Appl. Probab.}, 13\penalty0 (3):\penalty0 1119--1150, 2003.

\bibitem{har-tak}
J.~M. Harrison and M.~I. Taksar.
\newblock Instantaneous control of {B}rownian motion.
\newblock {\em Math. Oper. Res.}, 8\penalty0 (3):\penalty0 439--453, 1983.

\bibitem{Har-Van}
J.~M. Harrison and J.~A. Van~Mieghem.
\newblock Dynamic control of {B}rownian networks: state space collapse and
  equivalent workload formulations.
\newblock {\em Ann. Appl. Probab.}, 7\penalty0 (3):\penalty0 747--771, 1997.

\bibitem{Har-Wil}
J.~M. Harrison and R.~J. Williams.
\newblock Workload reduction of a generalized {B}rownian network.
\newblock {\em Ann. Appl. Probab.}, 15\penalty0 (4):\penalty0 2255--2295, 2005.

\bibitem{jac-shi}
J.~Jacod and A.~N. Shiryaev.
\newblock {\em Limit theorems for stochastic processes}, volume 288 of {\em
  Grundlehren der Mathematischen Wissenschaften [Fundamental Principles of
  Mathematical Sciences]}.
\newblock Springer-Verlag, Berlin, 1987.
\newblock ISBN 3-540-17882-1.
\newblock xviii+601 pp.

\bibitem{KLRS}
L.~Kruk, J.~Lehoczky, K.~Ramanan, and S.~Shreve.
\newblock An explicit formula for the {S}korokhod map on {$[0,a]$}.
\newblock {\em Ann. Probab.}, 35\penalty0 (5):\penalty0 1740--1768, 2007.

\bibitem{pw-little}
R.~Mazumdar, V.~Badrinath, F.~Guillemin, and C.~Rosenberg.
\newblock A note on the pathwise version of {L}ittle's formula.
\newblock {\em Oper. Res. Lett.}, 14\penalty0 (1):\penalty0 19--24, 1993.

\bibitem{PKH}
E.~Plambeck, S.~Kumar, and J.~M. Harrison.
\newblock A multiclass queue in heavy traffic with throughput time constraints:
  asymptotically optimal dynamic controls.
\newblock {\em Queueing Syst.}, 39\penalty0 (1):\penalty0 23--54, 2001.

\bibitem{plamb}
E.~L. Plambeck.
\newblock Optimal leadtime differentiation via diffusion approximations.
\newblock {\em Oper. Res.}, 52\penalty0 (2):\penalty0 213--228, 2004.

\bibitem{rei-snap}
M.~I. Reiman.
\newblock The heavy traffic diffusion approximation for sojourn times in
  {J}ackson networks.
\newblock In {\em Applied probability---computer science: the interface, {V}ol.
  {II} ({B}oca {R}aton, {F}la., 1981)}, volume~3 of {\em Progr. Comput. Sci.},
  pages 409--421. Birkh\"auser Boston, Boston, MA, 1982.

\bibitem{rub-ata}
M.~Rubino and B.~Ata.
\newblock Dynamic control of a make-to-order, parallel-server system with
  cancellations.
\newblock {\em Oper. Res.}, 57\penalty0 (1):\penalty0 94--108, 2009.

\bibitem{shi-ata-cid}
M.~Shifrin, R.~Atar, and I.~Cidon.
\newblock Optimal scheduling in the hybrid-cloud.
\newblock In {\em IFIP/IEEE International symposium on integrated network
  management}. IEEE, 2013, to appear.

\bibitem{smith}
W.~E. Smith.
\newblock Various optimizers for single-stage production.
\newblock {\em Naval Res. Logist. Quart.}, 3:\penalty0 59--66, 1956.

\bibitem{van}
J.~A. van Mieghem.
\newblock Dynamic scheduling with convex delay costs: the generalized {$c\mu$}
  rule.
\newblock {\em Ann. Appl. Probab.}, 5\penalty0 (3):\penalty0 809--833, 1995.

\bibitem{war-kum}
A.~R. Ward and S.~Kumar.
\newblock Asymptotically optimal admission control of a queue with impatient
  customers.
\newblock {\em Math. Oper. Res.}, 33\penalty0 (1):\penalty0 167--202, 2008.

\bibitem{whitt71}
W.~Whitt.
\newblock Weak convergence theorems for priority queues: {P}reemptive-resume
  discipline.
\newblock {\em J. Appl. Probability}, 8:\penalty0 74--94, 1971.

\bibitem{williams-ssc}
R.~J. Williams.
\newblock Diffusion approximations for open multiclass queueing networks:
  sufficient conditions involving state space collapse.
\newblock {\em Queueing Systems Theory Appl.}, 30\penalty0 (1-2):\penalty0
  27--88, 1998.

\end{thebibliography}

\end{document}